\documentclass{elsarticle}

\usepackage{hyperref}
\usepackage{bm}
\usepackage{float}
\usepackage{amsmath}
\usepackage{amsthm}
\usepackage{amsfonts}
\usepackage{graphicx}
\usepackage[usenames,dvipsnames]{xcolor}

\newtheorem{proposition}{Proposition}[section]
\newtheorem{minprinc}{Minimization Principle}[section]

\newcommand{\mcE}{\mathcal{E}}

\newcommand{\mcX}{\mathcal{X}}

\newcommand{\mcL}{\mathcal{L}}

\DeclareMathOperator*{\argmin}{arg\,min}

\newcommand{\Omg}{{\Omega}}
\newcommand{\Gam}{{\Gamma}}
\newcommand{\Gamu}{{\underline{\Gamma}}}
\newcommand{\Gami}{{\Gamma_1}}
\newcommand{\gam}{{\gamma}}

\newcommand{\ds}{\displaystyle}

\def \xb{\bm{x}}
\def \yb{\bm{y}}

\def \nb{\mathbf{n}}

\def\Omg{{\Omega}}

\usepackage{color}

\usepackage{setspace}
\begin{document}

\begin{frontmatter}

\title{An energy-based coupling approach to nonlocal interface problems}

\author{Giacomo Capodaglio \fnref{fsuFootnote}}

\author{Marta D'Elia \fnref{sandiaFootnote}}

\author{Pavel Bochev \fnref{sandiaFootnote}}

\author{Max Gunzburger \fnref{fsuFootnote}}

\fntext[fsuFootnote]{Department of Scientific Computing,
Florida State University,
e-mails: gcapodaglio, mgunzburger@fsu.edu}
\fntext[sandiaFootnote]{Center for Computing Research,
Sandia National Labs,
e-mail: pbboche, mdelia@sandia.gov}

\begin{abstract}
Nonlocal models provide accurate representations of physical phenomena ranging from fracture mechanics to complex subsurface flows, settings in which traditional partial differential equation models fail to capture effects caused by long-range forces at the microscale and mesoscale. However, the application of nonlocal models to problems involving interfaces, such as multimaterial simulations and fluid-structure interaction, is hampered by the lack of a physically consistent interface theory which is needed to support numerical developments and, among other features, reduces to classical models in the limit as the extent of nonlocal interactions vanish.
In this paper, we use an energy-based approach to develop a formulation of a nonlocal interface problem which provides a physically consistent extension of the classical perfect interface formulation for partial differential equations. Numerical examples in one and two dimensions validate the proposed framework and demonstrate the scope of our theory.
\end{abstract}
 

\end{frontmatter}

\section{Introduction}

Nonlocal models can accurately describe physical phenomena arising from long-range forces at the microscale and mesoscale. Such phenomena cannot be accounted for by partial differential equation (PDE) models in which the interaction is limited to points that are in direct contact with each other. Nonlocal models are represented mathematically by integral operators which are better suited to capture interactions occurring across a distance. In particular, physically consistent nonlocal models have been defined that allow for the treatment of nonstandard effects such discontinuous solutions and, more generally, can capture nonstandard effects such as multiscale behaviors of the displacement in continuum mechanics and anomalous diffusion processes in, e.g., subsurface flow applications. In fact, specific examples in which long-range forces are essential for accurate
predictive simulations can be found in a diverse spectrum of scientific applications such as anomalous subsurface transport \cite{benson2000application, schumer2001eulerian, schumer2003multiscaling, delgoshaie2015non}, fracture mechanics \cite{silling2000reformulation, ha2011characteristics, littlewood2010simulation}, image processing \cite{buades2010image, gilboa2007nonlocal, gilboa2008nonlocal, lou2010image}, magnetohydrodynamics \cite{schekochihin2008mhd},  multiscale and multiphysics systems \cite{askari2008peridynamics, alali2012multiscale}, phase transitions \cite{bates1999integrodifferential, fife2003some}, and stochastic processes \cite{meerschaert2011stochastic, d2017nonlocal, burch2014exit}.

 Although  research efforts devoted to nonlocal models have recently intensified, there are few such examples specifically devoted to nonlocal interface (NLI) problems, with \cite{alali2015peridynamics, seleson2013interface} being perhaps the only published works in this direction. 
In contrast, the extant literature on interface problems governed by local PDEs is vast and ranges from standard transmission conditions \cite{Hansbo_02_CMAME}  to problems with more complicated conditions such as the Beavers--Joseph--Saffman--Jones approximation \cite{Chen_11_SINUM}, or the Beavers-Joseph conditions \cite{Cao_10_CMS} in Stokes-Darcy models. A thorough review of this literature is beyond the scope of the paper and we limit ourselves to just these references to provide some illustrative examples.

 However, the NLI formulations in these papers do not demonstrate convergence, as the extent of nonlocal interactions vanish, of their nonlocal models and their solutions to their classical local PDE counterparts; such convergence is a litmus test that nonlocal models must pass in many settings, including the ones mentioned above. Furthermore, these papers do not establish, for the nonlocal models they consider, a mathematically rigorous well-posedness theory nor do they provide rigorous error estimation for approximate solutions. The absence of such a mathematical framework has hampered the wider adoption of nonlocal models in science and engineering applications. In this paper we also do not provide such theories but will do so in a follow-up paper.

Here, our goals are first to define a new NLI model. The NLI model we introduce is motivated by an energy-based description of classical local interface problems in which the local energy of the system is minimized subject to constraints modeling the physics occurring across the local interface. The next goal is to show that the new model does indeed pass the litmus test mentioned above. We then also provide results of numerical experiments which illustrate the behavior of solutions with respect to changes in grid size and extent of nonlocal interactions and which also illustrate the convergence of nonlocal models and their solution to their local counterparts. Being the first such efforts in these directions, we develop the NLI model in the simple context of the nonlocal counterpart of the local interface problem for the Poisson equation. Doing so allows for a clear exposition of the important features of the NLI model and of its solution and of comparisons with local interface models. Comments regarding generalization to more complex settings are made at the end of the paper.

The paper is structured as follows. In Section \ref{sec:loc_int_prob} we review the classical energy-based formulation of a local interface problem for second-order elliptic PDEs; this review provides the template for the development of an energy-based formulation of a NLI problem. In Section \ref{sec:non_vol_cons}, we review some fundamental notions about nonlocal problems. The two reviews given in Sections \ref{sec:loc_int_prob} and 
\ref{sec:non_vol_cons} set us up for introducing, in Section \ref{sec:geometry-energy-kernel}, the NLI model we consider in the rest of the paper. Then, in Section \ref{loc_lims}, we prove that that model passes the litmus test and, in Section \ref{sec:num_res}, we present our numerical results. Finally, in Section \ref{sec:conclusions}, concluding remarks are provided, including prospects for future work.

\section{A local, PDE-based, interface problem}\label{sec:loc_int_prob}

In this section, we review a classical energy-based formulation of local interface problems for second-order elliptic PDEs.

\begin{figure}[h!]
\centering
\includegraphics[width=0.25\textwidth]{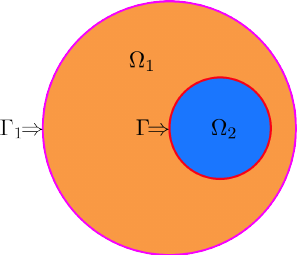}
\hspace{1cm}
\includegraphics[width=0.25\textwidth]{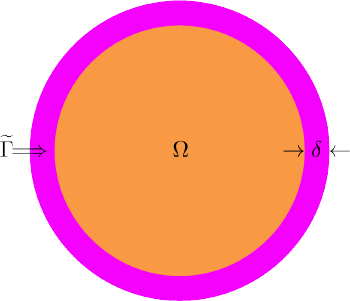}
\caption{Left: illustration of the geometric configuration for the local interface problem. Right: illustration of the geometric configuration for the nonlocal volume-constrained problem.}
\label{fig:loc_int}
\end{figure}

Let $\Omega_1$ and $\Omega_2$ be two disjoint open and bounded subsets of $\mathbb{R}^n$, $n=2,3$, with $\Omega_2$ surrounded by $\Omega_1$ as illustrated in Figure \ref{fig:loc_int}-left.\footnote{Note that in one dimension, this description is not possible. A one-dimensional configuration is described Section \ref{sec:num_res} and used for several numerical tests. Note also that, for practical implementation reasons, the two-dimensional domain configuration used for the numerical tests is different and consists of rectangular domains, see Figure \ref{2D}, top left.} The interface, i.e., the common boundary between the domains, is defined as $\Gamma=\overline\Omega_1\cap\overline\Omega_2$. With $\partial\Omega_i$ denoting the boundary of $\Omega_i$, we let $\Gamma_i=\partial\Omega_i\setminus\Gamma$. Note that $\partial\Omega_2=\Gamma$. The geometric entities so introduced are illustrated in Figure \ref{fig:loc_int}-left.

\subsection{Local energy minimization principle}\label{sec:LEMP}
Consider the following (local) energy functional
\begin{equation}\label{eq:energy_loc_int}
\begin{aligned}
    {\mathcal E}(u_1,u_2;f_1,f_2)
=&\frac12 \int_{\Omega_1} \kappa_1(\xb)|\nabla u_1(\xb)|^2d\xb -\int_{\Omega_1} f_1(\xb)u_1(\xb)d\xb
\\&\qquad+\frac12\int_{\Omega_2} \kappa_2(\xb)|\nabla u_2(\xb)|^2d\xb  - \int_{\Omega_2} f_2(\xb)u_2(\xb)d\xb,
\end{aligned}
\end{equation}
where the functions $\kappa_1$ and $\kappa_2$ represent the different material properties of the two domains and are assumed to be positive and bounded from below. The functions $f_1$, $f_2$ are known.

We obtain a particular instance of a local interface problem by choosing a specific constrained minimization setting for \eqref{eq:energy_loc_int}. To this end, let us define the following energy spaces, for $i=1,2$ 
\begin{equation}
\begin{aligned}
W_i &= \{w\in H^1(\Omega_i) \,:\, |||w|||_i < \infty\}, \quad \mbox{where} \quad |||w|||^2_i=\|\nabla w\|^2_{0,\Omega_i}+\| w\|^2_{0,\Omega_i}, \\
W^c_i &= \{w\in W_i \,:\, w(\xb)=0 \,\, \mbox{on} \,\, \Gamma_i\}.
\end{aligned}
\end{equation}
Note that, due to the configuration of the domains depicted in Figure \ref{fig:loc_int}, left, $W^c_2=W_2$. Tensor product spaces are then defined as $W=W_1 \times W_2$ and $W^c=W^c_1 \times W^c_2$.
\begin{minprinc} \label{eq:local-min-pr}
Given $g_1$, $f_1$, $f_2$, $\kappa_1$ and $\kappa_2$, find $(u_1,u_2)\in W$ such that
$$
    (u_1,u_2) = \argmin\limits_{(v_1,v_2)\in W} {\mathcal E}(v_1,v_2;f_1,f_2), 
$$
subject to the constraints
\begin{equation}\label{eq:cond_en_princ}
\begin{cases}
u_1(\xb)= g_1(\xb)  &\quad   \xb\in\Gamma_1,
\\
u_1(\xb) =  u_2(\xb) &\quad   \xb\in\Gamma.
\end{cases}
\end{equation}
\end{minprinc}
\noindent Here, $g_1$ is a known Dirichlet data function defined on $\Gamma_1$. The second constraint in \eqref{eq:cond_en_princ}, i.e., the continuity of the states across the interface, is a modeling assumption about the physics of the interface, which gives rise to a specific flavor of a local interface problem. 

\subsection{Weak formulation}
The Euler-Lagrange equation corresponding to the Minimization Principle \ref{eq:local-min-pr} is given by the following weak variational equation:
 find $(u_1,u_2)\in W$ satisfying the constraints in \eqref{eq:cond_en_princ} and such that
\begin{equation}\label{eq:weak_recovered}
\begin{aligned}
&\qquad\int_{\Omega_1} \kappa_1(\xb)\nabla u_1(\xb)\cdot \nabla v_1(\xb)d\xb  
+\int_{\Omega_2} \kappa_2(\xb)\nabla u_2(\xb)\cdot \nabla v_2(\xb)d\xb  
=\int_{\Omega_1} f_1(\xb)v_1(\xb)d\xb 
+\int_{\Omega_2} f_2(\xb)v_2(\xb)d\xb 
\end{aligned}
\end{equation}
for all $(v_1,v_2)\in W^c$ satisfying $v_1(\xb) = v_2(\xb)$ on $\Gamma$.

\subsection{Strong formulation}
As usual, we derive the strong form of the interface problem from the weak formulation \eqref{eq:weak_recovered} by assuming that $u_1$ and $u_2$ are sufficiently regular. Collecting terms, integrating by parts, and taking into account that $v_1\in W^c$ yields
\begin{align}\label{eq:temp}
   \int_{\Omega_1}& \Big(-\nabla\cdot(\kappa_1(\xb)\nabla u_1(\xb))  - f_1(\xb)\Big) v_1(\xb)  d\xb
  + \int_{\Omega_2} \Big(-\nabla\cdot (\kappa_2(\xb)\nabla u_2(\xb)) - f_2(\xb) \Big) v_2(\xb)d\xb \\ \nonumber
&+\int_{\Gamma} v_1(\xb) \kappa_1(\xb)\nabla u_1(\xb)\cdot \nb_1 d\xb
  + \int_{\Gamma} v_2(\xb) \kappa_2(\xb)\nabla u_2(\xb)\cdot \nb_2 d\xb
= 0,
\end{align}
where $\nb_1$ and $\nb_2$ are unit vectors, normal to the interface $\Gamma$, pointing outward from $\Omega_1$ and $\Omega_2$, respectively.
Because $v_1$ and $v_2$ are arbitrary on $\Omega_1\cup\Gamma$ and $\Omega_2\cup\Gamma$, respectively, we may first set $v_1$ arbitrary on $\Omega_1$, $v_1=0$ on $\Gamma$, and $v_2=0$ on $\Omega_2\cup\Gamma$ and then set $v_2$ arbitrary on $\Omega_2$, $v_2=0$ on $\Gamma$, and $v_1=0$ on $\Omega_1\cup\Gamma$ to obtain from \eqref{eq:temp} the strong forms $-\nabla\cdot(\kappa_i(\xb)\nabla u_i(\xb))  = f_i(\xb)$ on $\Omega_i$, $i=1,2$, of the subdomain equations. Substituting these equations back into \eqref{eq:temp} leaves us with 
$$
\int_{\Gamma} v_1(\xb) \kappa_1(\xb)\nabla u_1(\xb)\cdot \nb_1 d\xb + \int_{\Gamma} v_2(\xb) \kappa_2(\xb)\nabla u_2(\xb)\cdot \nb_2 d\xb = 0.
$$
Using that $v_1=v_2$ on $\Gamma$ we then recover from this equation the \emph{flux continuity condition} 
$\kappa_1(\xb)\nabla u_1(\xb)\cdot{\bf n}_1 + \kappa_2(\xb)\nabla u_2(\xb)\cdot{\bf n}_2 = 0$.
Thus, the strong form of the local interface problem corresponding to the Minimization Principle \eqref{eq:local-min-pr} is given by
\begin{align}\label{eq:strong_loc_int}
-\nabla\cdot \big(\kappa_1(\xb)\nabla u_1(\xb)\big) &= f_1(\xb) &
\xb\in{\Omega_1}, \\
u_1(\xb)&=g_1(\xb)  & \xb\in\Gamma_1, \label{eq:strong_loc_int-BC1}\\
-\nabla\cdot \big(\kappa_2(\xb)\nabla u_2(\xb)\big) &= f_2(\xb) & 
\xb\in\Omega_2,\\
u_1(\xb) &=  u_2(\xb) & \xb\in\Gamma, \label{eq:strong_loc_int-state-cont}\\
\kappa_1(\xb)\nabla u_1(\xb)\cdot{\bf n}_1 + 
\kappa_2(\xb)\nabla u_2(\xb)\cdot{\bf n}_2 &= 0 & \xb\in\Gamma.
\label{eq:strong_loc_int-flux-cont}
\end{align}
We note that the strong form of the interface problem contains the flux continuity condition \eqref{eq:strong_loc_int-flux-cont} that was not explicitly present in Minimization Principle \ref{eq:energy_loc_int}. This condition is a consequence of the specific form of the energy functional \eqref{eq:energy_loc_int} and of the constraint \eqref{eq:strong_loc_int-state-cont}, both of which are modeling assumptions about the physics local interface problem. 
Interfaces for which both the jumps in the state and in the normal flux are zero across the interface are known as \emph{perfect interfaces} \cite{Javili_14_CMAME}. In contrast, interfaces for which one or both of these quantities are discontinuous across the interface are known as \emph{imperfect interfaces}; see, e.g. \cite{Javili_14_CMAME}. 

One practically important example are interfaces in which  the flux is driven by a jump in the primal variable. Such imperfect interfaces are employed across vastly different scales to model, e.g., interfacial thermal resistance at the nanoscale in semiconductor devices  \cite{Ming_08_JAP,Berber_00_PRL}, as well as  the flux exchange between the ocean and the atmosphere \cite{Lemarie_15_PCS} in climate models. In both cases the interface physics is not fully understood and is modeled by constitutive ``closure'' relations. Furthermore, It is worth pointing out that mathematical and numerical analysis of such conditions is somewhat limited even in the local setting; see, e.g., \cite{Lions_95_JMPA} and \cite{Kwak_17_IJNME} for some of the available examples. 
For this reason, we leave extension of imperfect interfaces to nonlocal problems for future work and focus instead solely on nonlocal generalizations of perfect interfaces.

\section{Nonlocal volume constrained problems}\label{sec:non_vol_cons}

In this section, we review the fundamentals for a nonlocal problem which is the nonlocal counterpart of the local Poisson PDE. Let $\Omega$ be an open and bounded subset of $\mathbb{R}^n$. Given a positive real number $\delta$, often referred to as the {\em horizon} or {\em interaction radius}, we define the {\it interaction domain} $\widetilde{\Gamma}$ associated with $\Omega$ as 
\begin{align}
    \widetilde{\Gamma} = \{\yb \in  \mathbb{R}^n  \setminus \Omega\,:\, |\xb - \yb| \leq \delta \,, \mbox{for} \, \xb \in \Omega\}.
\end{align}
Note that $\widetilde{\Gamma}$, as all other entities defined in this section, depend on $\delta$, so that, for the sake of economy of notation, that dependence is not explicitly included in those entities.
Figure \ref{fig:loc_int}-right illustrates an example of a two-dimensional domain $\Omega$ and its interaction domain $\widetilde{\Gamma}$ having thickness $\delta$, i.e., having non-zero volume in ${\mathbb R}^n$.

\subsection{Nonlocal energy minimization principle}

We use an energy-based characterization of nonlocal volume constrained problems which mirrors the Dirichlet principle for the classical gradient operator. Specifically, we seek the states of the nonlocal model as suitably constrained minimizers of the nonlocal energy functional
\begin{align}\label{eq:energy_nonloc}
     \mcE(u;\gamma,f)=\frac12\int_{\Omega\cup\widetilde\Gamma}  \int_{\Omega\cup\widetilde\Gamma}\big|u(\yb) - u(\xb)\big|^2\gamma(\xb,\yb)\, d\yb  d\xb
- \int_\Omega f(\xb)u(\xb)d\xb. 
 \end{align}
The function $\gamma$ is referred to as the {\it{kernel}} and is required to satisfy
\begin{align}
     \gamma(\xb,\yb) = \gamma(\yb,\xb), \quad \mbox{for} \, \xb,\yb \in \Omega \cup \widetilde{\Gamma}.
 \end{align}
Let us define the following function spaces
\begin{equation}
 \begin{aligned}
&V = \{v\in L^2(\Omega\cup\widetilde\Gamma) \,\,:\,\, |||v|||<\infty\},\\
\quad &\mbox{where} \quad  |||v|||^2 =  \int_{\Omega\cup\widetilde\Gamma}  \int_{\Omega\cup\widetilde\Gamma}|v(\yb) - v(\xb)|^2\gamma(\xb,\yb)\,d\yb  d\xb + \|v\|^2_{L^2(\Omega\cup\widetilde\Gamma)},\\
&V^c = \{v\in V\,\,:\,\, v=0 \,\,\mbox{for $\xb\in\widetilde\Gamma$}\}.
\end{aligned}
\end{equation}
\begin{minprinc} \label{eq:min-prin-nl}
Given $\gamma$, $f$, and $g$, find $u \in V$ such that
$$
u = \argmin\limits_{v\in V} {\mathcal E}(v;\gamma,f), 
$$
subject to $u(\xb)=g(\xb)$ on $\widetilde{\Gamma}$.
\end{minprinc}
\noindent Note that the constraint $u(\xb)\!=\!g(\xb)$ is applied on the interaction domain $\widetilde\Gamma$ that has nonzero volume, in contrast to the local case for which such a Dirichlet constraint is applied on the boundary $\partial\Omega$. We refer to the constraint in \eqref{eq:min-prin-nl} as a Dirichlet {\em volume constraint}.
%

\subsection{Weak formulation}
The necessary optimality condition of the Minimization Principle \ref{eq:min-prin-nl} is given by the following variational equation: find $u \in V$ such that $u(\xb)\!=\!g(\xb)$ on $\widetilde{\Gamma}$ and
\begin{equation}\label{eq:weak_recovered_2}
\begin{aligned}\int_{\Omega\cup\widetilde\Gamma}  \int_{\Omega\cup\widetilde\Gamma}\big(v(\yb) - v(\xb)\big)\big(u(\yb) - u(\xb)\big)\gamma(\xb,\yb)d\yb  d\xb
=  \int_\Omega v(\xb)f(\xb)d\xb \qquad \forall \,v \in V^c.
\end{aligned}
\end{equation}

\subsection{Strong formulation}
To state the strong form of \eqref{eq:weak_recovered_2} we recall the nonlocal diffusion operator
 \begin{align}\label{eq:NL-diff}
    \mcL u(\xb) = 2\int_{\Omega\cup\widetilde\Gamma} \big(u(\yb) - u(\xb)\big)\gamma(\xb,\yb)d\yb.
 \end{align}
and the nonlocal Green's identity \cite{Du_13_MMMAS}
\begin{equation}
 \begin{aligned}\label{eq:Green_nonloc}
       \int_\Omega  v(\xb) \mcL u(\xb) d\xb = &-
  \int_{\Omega\cup\widetilde\Gamma}  \int_{\Omega\cup\widetilde\Gamma}\big(v(\yb) - v(\xb)\big)\big(u(\yb) - u(\xb)\big)\gamma(\xb,\yb)d\yb  d\xb \\
  & -2\int_{\widetilde\Gamma}  \int_{\Omega\cup\widetilde\Gamma} v(\xb)\big(u(\yb) - u(\xb)\big)\gamma(\xb,\yb)d\yb  d\xb.
 \end{aligned}
 \end{equation}
Using \eqref{eq:Green_nonloc} and the fact that  $v=0$ on $\widetilde{\Gamma}$ one can transform \eqref{eq:weak_recovered_2} into the following equation
  \begin{align}\label{eq:weak_1_nonloc}
     \int_\Omega v(\xb) \mcL u (\xb) d\xb +  \int_\Omega v(\xb) f(\xb) d\xb = 0.
 \end{align}
Since $v$ is arbitrary on $\Omega$ one then easily obtains the strong form of the nonlocal volume constrained problem
\begin{equation}\label{eq:strong_nonloc}
    \begin{aligned}
    \displaystyle -\mcL u(\xb) &= f(\xb) & \xb\in\Omega
\\[.5ex]
u(\xb) &= g(\xb) &  \xb\in\widetilde\Gamma 
    \end{aligned}
\end{equation}
which is the nonlocal counterpart of the Dirichlet problem for second-order elliptic PDEs.


\section{Nonlocal interface problem} \label{sec:geometry-energy-kernel}

In this section we define a NLI model that provides a physically consistent extension of the perfect local interface problem of Section \ref{sec:loc_int_prob}.
%

\begin{figure}[t!]
\centering
\includegraphics[width=0.24\textwidth]{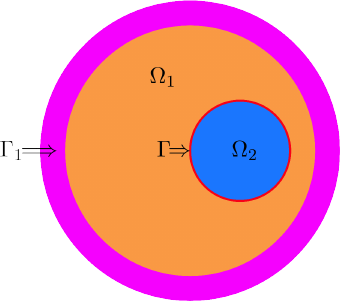}
\includegraphics[width=0.3\textwidth]{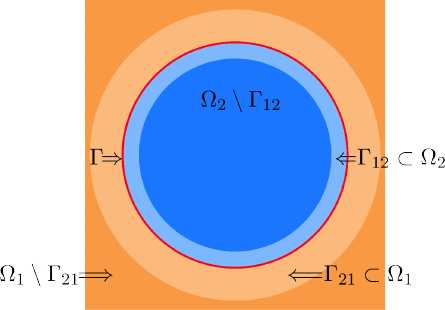}
\includegraphics[width=0.3\textwidth]{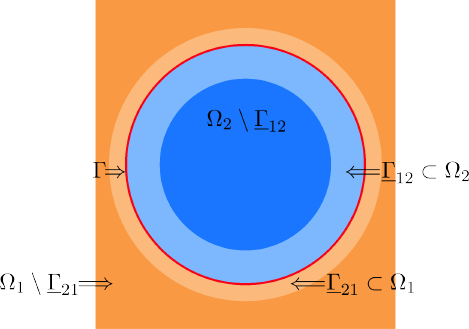}
\caption{Left: geometric configuration for the nonlocal interface problem. Center: illustration of the subdomains $\Gamma_{12}$ and $\Gamma_{21}$. Right: illustration of the subdomains $\underline\Gamma_{12}$ and $\underline\Gamma_{21}$.}
\label{fig:fig2}
\end{figure}

%
In the NLI model, the domains $\Omega_i$, $i=1,2$, interact with each other through regions that have nonzero volume. These regions are defined in different ways depending on the relative location of $\xb$ and $\yb$ and on two horizon parameters $\delta_1$ and $\delta_2$. Specifically, for $\delta_1 >0$ and $\delta_2 > 0$, we let
\begin{itemize}
\item $\Gami = \{\yb \in 
\mathbb{R}^n \setminus \Omega_1 \,:\, |\xb - \yb| \leq \delta_1 \,, \mbox{for} \, \xb \in \Omega_1\}$ : the (external) interaction domain of $\Omg_1$.
\item $\Gam_{12} = \{\yb \in \Omega_2 \,:\, |\xb - \yb| \leq \delta_1 \,, \mbox{for} \, \xb \in \Omega_1\}$: the subregion of $\Omg_2$ that interacts with $\Omg_1$ when $\xb\in\Omg_1$. 
\item $\Gam_{21}= \{\yb \in \Omega_1  \,:\, |\xb - \yb| \leq \delta_2 \,, \mbox{for} \, \xb \in \Omega_2\}$: the subregion of $\Omg_1$ that interacts with $\Omg_2$ when $\xb\in\Omg_2$.
\item $\Gamu_{12}= \{\yb \in \Omega_2 \,:\, |\xb - \yb| \leq \delta_2 \,, \mbox{for} \, \xb \in \Omega_1\}$: the analogue of $\Gam_{21}$ on $\Omega_2$.
\item $\Gamu_{21}= \{\yb \in \Omega_1\,:\, |\xb - \yb| \leq \delta_1 \,, \mbox{for} \, \xb \in \Omega_2\}$: the analogue of $\Gam_{12}$ on $\Omega_1$.
\end{itemize}
Figure \ref{fig:fig2} illustrates the geometric configuration for the nonlocal interface problem and the various subdomains involved.
We also introduce the set $\Gamma^* = \Gamma_{ij} \cup \Gamu_{ji}$ if $\delta_i \geq \delta_j$. Note that by definition $\Gam_{ij} \subset \Gam^*$ and $\Gamu_{ij} \subset \Gam^*$ for any $i$ and $j$.
The kernel function $\gamma$ is defined for all $(\xb,\yb) \in (\Omega \cup \Gamma_1) \times (\Omega \cup \Gamma_1)$ as follows:
\begin{equation}\label{eq:kernel}
\gam(\xb,\yb)=\left\{\begin{array}{ll}
\gam_{11}(\xb,\yb) = C_{11}(\delta_1) \mcX_{B_1(\xb)}(\yb) 
                   & \xb\in \Omega_1\cup\Gami, \; \yb\in\Omega_1\cup\Gami\\[2mm]
\gam_{12}(\xb,\yb) = C_{12}(\delta_1) \mcX_{B_1(\xb)}(\yb) 
                   & \xb\in\Omega_1\cup\Gami, \; \yb\in\Omg_2\\[2mm]
\gam_{21}(\xb,\yb) = C_{21}(\delta_2) \mcX_{B_2(\xb)}(\yb) 
                   & \xb\in\Omg_2, \; \yb\in\Omega_1 \cup\Gami\\[2mm]
\gam_{22}(\xb,\yb) = C_{22}(\delta_2) \mcX_{B_2(\xb)}(\yb) 
                   & \xb\in\Omg_2, \; \yb\in\Omg_2,                   
\end{array}\right.
\end{equation}
where $\mcX_S$ is the indicator function on the set $S$ and $B_i(\xb)$ is the ball of radius $\delta_i$ centered at $\xb$. Specific choices of $C_{ij}(\delta_i)$, $i,j=1,2$, are discussed in Section \ref{loc_lims} and Section \ref{sec:num_res}.

\subsection{Nonlocal energy minimization principle}
Mirroring the energy-based description of the local perfect interface problem we start with defining the nonlocal energy of the system as follows
\begin{equation}\label{eq:energy-unsplit}
\mathcal{E}(u;\gamma, f)=\ds\frac12\int\limits_{\Omg\cup\Gami}\int\limits_{\Omg\cup\Gami}
(u(\xb)-u(\yb))^2 \gam(\xb,\yb)\,d\yb\,d\xb -
\int\limits_{\Omg} f(\xb) u(\xb)\, d\xb
\end{equation}
Using that $\Omega_1$ and $\Omega_2$ are disjoint, we split the energy in two parts, associating the first to $u_1$ and the second to $u_2$, i.e.
\begin{equation}\label{eq:energy-split}
\begin{aligned}
\mathcal{E}_s(u_1,u_2;\gamma, f)
& = \ds\frac12\int\limits_{\Omega_1\cup\Gami}\int\limits_{\Omg\cup\Gami}
    (u_1(\xb)-u_1(\yb))^2 \gam(\xb,\yb)\,d\yb\,d\xb - \int\limits_{\Omega_1} f u_1\, d\xb\\[2mm]
& + \frac12\int\limits_{\Omg_2}\int\limits_{\Omg\cup\Gami}
    (u_2(\xb)-u_2(\yb))^2 \gam(\xb,\yb)\,d\yb\,d\xb - \int\limits_{\Omg_2} f u_2\, d\xb.
\end{aligned}
\end{equation}
Let us define the following function spaces
 \begin{equation}
 \begin{aligned}
        &W_1 = \{ w \in L^2(\Omega_1\cup\Gamma_1)\,\,\,:\,\,\,
   |||w|||_1^2 + \|w\|^2_{L^2(\Omega_1\cup\Gamma_1)}
   <\infty     \big\}, \\
        \quad &\mbox{where} \quad |||w|||_1^2 = \frac12\int_{\Omega_1\cup\Gamma_1}\int_{\Omega_1\cup\Gamma_1}
|u_1(\yb)-u_1(\xb)|^2\gamma_{11}(\xb,\yb)d\yb d\xb,\\
        &W_2 = \{ w \in L^2(\Omega_2)\,\,\,:\,\,\,
   |||w|||_2^2 + \|w\|^2_{L^2(\Omg_2)}
   <\infty     \big\}, \\
        \quad &\mbox{where} \quad |||w|||_2^2 = \frac12\int_{\Omega_2}\int_{\Omega_2}
|u_2(\yb)-u_2(\xb)|^2\gamma_{22}(\xb,\yb)d\yb d\xb.\\
 \end{aligned}
 \end{equation}
 The constrained space $W^c_1$ is defined as
 \begin{align}
       W_1^c = \big\{  w\in W_1\,\,\,:\,\,\, \mbox{$w(\xb) = 0$ on $\Gamma_1$}\}.
 \end{align}
We also introduce the tensor product spaces $W=W_1 \times W_2$ and $W^c = W^c_1 \times W_2$. Note that we used the same notation as in the local case even though these are different functional spaces.
\begin{minprinc} \label{eq:minprin-NLI}
Given $\gamma$, $g_1$ and $f$, find $(u_1,u_2)\in W$ such that
$$
    (u_1,u_2) = \argmin\limits_{(v_1,v_2)\in W} {\mathcal E}(v_1,v_2;\gamma,f), 
$$
subject to the constraints
\begin{align}\label{eq:cond_en_princ_nl}
\begin{cases}
u_1(\xb)= g_1(\xb)  &\quad   \xb\in\Gamma_1,
\\
u_1(\xb) =  u_2(\xb) &\quad   \xb\in \Gamma^*.
\end{cases}
\end{align}
\end{minprinc}
The second constraint in \eqref{eq:cond_en_princ_nl} can be viewed as a generalization of the state continuity constraint in  \eqref{eq:cond_en_princ} much like the volume constraint in \eqref{eq:strong_nonloc} generalizes the standard Dirichlet boundary condition.

\subsection{Weak formulation}

From the Euler-Lagrange equation corresponding to the Minimization Principle \ref{eq:minprin-NLI}, we derive a weak variational equation for the NLI model given in the following proposition. The derivation is given in \ref{sec:weak-form-proof}.

\begin{proposition}\label{lem:weak-form}
A weak formulation corresponding to the Minimization Principle \ref{eq:minprin-NLI} is given by the following weak variational equation: find $u_1\in W_1$ and  $u_2\in W_2$ such that \eqref{eq:cond_en_princ_nl} and $v_1(\xb) =  v_2(\xb)$ for $\xb\in \Gamma^*$ are satisfied and
\begin{equation}\label{eq:weak-omega1}
\begin{aligned}
 &\int\limits_{\Omega_1 \cup\Gami}\int\limits_{\Omega_1\cup\Gami}
    (u_1(\xb)-u_1(\yb))(v_1(\xb)-v_1(\yb)) \gamma_{11}(\xb,\yb)d\yb\,d\xb
+ \int\limits_{\Omega_1}\int\limits_{\Omg_2}
(u_1(\xb)-u_2(\yb))v_1(\xb) \gamma_{12}(\xb,\yb)d\yb\,d\xb
\\
&\qquad-\int\limits_{\Omg_2}\int\limits_{\Omega_1} (u_2(\xb)-u_1(\yb))v_1(\yb) \gamma_{21}(\xb,\yb)d\yb\,d\xb
= \int\limits_{\Omega_1} f(\xb) v_1(\xb)\, d\xb 
\qquad\mbox{for all $v_1(\xb)\in W_1^c$}
\end{aligned}
\end{equation}
and
\begin{equation}\label{eq:weak-omega2}
\begin{aligned}
&\int\limits_{\Omg_2}\int\limits_{\Omg_2} (u_2(\xb)-u_2(\yb))(v_2(\xb)-v_2(\yb)) \gamma_{22}(\xb,\yb)d\yb\,d\xb
+\int\limits_{\Omg_2}\int\limits_{\Omega_1} (u_1(\yb)-u_2(\xb)v_2(\xb)) \gamma_{21}(\xb,\yb)d\yb\,d\xb
\\&\qquad- \int\limits_{\Omega_1}\int\limits_{\Omg_2}
(u_1(\xb)-u_2(\yb))v_2(\yb) \gamma_{12}(\xb,\yb)d\yb\,d\xb
=\int\limits_{\Omg_2} f(\xb) v_2(\xb)\, d\xb
\qquad\mbox{for all $v_2(\xb)\in W_2$.}
\end{aligned}
\end{equation}
\end{proposition}

\subsection{Strong formulation}

The strong formulation of the NLI model can be derived from the weak formulation as shown in the following proposition whose proof is given in \ref{sec:strong-form-proof}.
\begin{proposition}\label{lem:strong-form}
The strong form of the NLI model associated with the weak formulation \eqref{eq:weak-omega1} and \eqref{eq:weak-omega2} is given by the nonlocal subdomain equations
\begin{equation}\label{eq:strong1}
-2\!\!\!\int\limits_{\Omega_1\cup\Gam_1} (u_1(\yb)-u_1(\xb))\gam_{11}(\xb,\yb)\,d\yb
-  \!\int\limits_{\Omg_2} (u_2(\yb)-u_1(\xb))\big(\gam_{12}(\xb,\yb)+\gam_{21}(\xb,\yb)\big)\,d\yb= f(\xb),
\end{equation}
for $\xb\in\Omega_1$, and
\begin{equation}\label{eq:strong2}
-2\!\int\limits_{\Omg_2} (u_2(\yb)-u_2(\xb))\gam_{22}(\xb,\yb)\,d\yb
-\!\int\limits_{\Omg_1} (u_1(\yb)-u_2(\xb))\big(\gam_{21}(\xb,\yb)+\gamma_{12}(\xb,\yb)\big)\,d\yb= f(\xb),
\end{equation}
for $\xb\in\Omg_2$. 
\end{proposition}

Isolating terms in \eqref{eq:strong1} and \eqref{eq:strong2} that do not interact with $\Omega_2$ and $\Omega_1$, respectively, we obtain the equations
\begin{equation}\label{eq:strong-domain}
\begin{array}{ll}
-2C_{11}(\delta_1)\ds\int_{B_1(\xb)} (u_1(\yb)-u_1(\xb))\,d\yb = f(\xb)
   & \quad\mbox{for $\xb\in\Omega_1\setminus\Gamu_{21}$}, \\[6mm] 
-2C_{22}(\delta_2)\ds\int_{B_2(\xb)} (u_2(\yb)-u_2(\xb))\,d\yb = f(\xb)
   &  \quad\mbox{for $\xb\in\Omg_2\setminus\Gamu_{12}$}, 
\end{array}
\end{equation}
{where we have used \eqref{eq:kernel} to substitute for $\gamma_{ii}(\xb,\yb)$ and also that, by construction, $(\Omega_1\cup\Gam_1)\cap B_1(\xb)=B_1(\xb)$ for $\xb\in\Omega_1\setminus\Gamu_{21}$ and $\Omega_2\cap B_2(\xb)=B_2(\xb)$ for $\xb\in\Omega_2\setminus\Gamu_{12}$.}
Collecting the remaining terms then yields 
\begin{equation}\label{eq:strong-interface1}
-2C_{11}(\delta_1)\!\!\!\!\!\!\!\!\int\limits_{(\Omega_1\cup\Gam_1)\cap B_1(\xb)}
\!\!\!\!\!\!\!\!\!\!
(u_1(\yb)-u_1(\xb))\,d\yb
-C_{12}(\delta_1)\!\!\!\!\!{\int\limits_{\Gamma_{12}\cap B_1(\xb)}}\!\!\!\!\!\! (u_2(\yb)-u_1(\xb))\,d\yb
-C_{21}(\delta_2)\!\!\!\!\!{\int\limits_{\underline\Gamma_{12}\cap B_2(\xb)}}\!\!\!\!\!\! (u_2(\yb)-u_1(\xb))\,d\yb= f(\xb) 
\end{equation}
for $\xb\in \Gamu_{21}$, and
\begin{equation}\label{eq:strong-interface2}
-2C_{22}(\delta_2)\!\!\!\!\int\limits_{\Omg_2\cap B_2(\xb)}\!\!\!\!\!\! (u_2(\yb)-u_2(\xb)) \,d\yb
-C_{21}(\delta_2)\!\!\!\!\!{\int\limits_{\Gamma_{21}\cap B_2(\xb)}}\!\!\!\!\!\!
(u_1(\yb)-u_2(\xb))\,d\yb  
-C_{12}(\delta_1)\!\!\!\!\!{\int\limits_{\underline\Gamma_{21}\cap B_1(\xb)}}\!\!\!\!\!\!
(u_1(\yb)-u_2(\xb))\,d\yb  =f(\xb)
\end{equation}
for $\xb\in \Gamu_{12}$, { where we have again used \eqref{eq:kernel}.
}

To summarize, the strong form of the nonlocal interface problem is comprised of the subdomain equations \eqref{eq:strong-domain}, the nonlocal flux interface conditions \eqref{eq:strong-interface1}  and \eqref{eq:strong-interface2}, and the volumetric constraints in \eqref{eq:cond_en_princ_nl}. Together, these nonlocal equations and interface conditions are the nonlocal analogue of the perfect interface PDE formulation \eqref{eq:strong_loc_int}--\eqref{eq:strong_loc_int-flux-cont}.


\section{Local limits of the nonlocal interface model}\label{loc_lims}
In the next proposition whose proof is given in  \ref{sec:local-limit-proof}, we show that the weak form of the NLI problem converges to the weak form of the local interface problem as the extent of the nonlocal interactions vanish, i.e. as $\delta_1$ and $\delta_2$ approach zero. This exercise results in unique choices for the constants $C_{ii}$, $i,j=1,2$, whereas it allows for some freedom in the choices of $C_{ij}$ for $i \neq j$. For the sake of clarity, we present the analysis in the two-dimensional setting; extension to the three-dimensional case is straightforward but involves significantly more cumbersome notation.

\begin{proposition}\label{lem:local-limit}
For $i=1,2$, and any positive constants $\widetilde C_{ij}$, provided $u_i(\xb)$ is sufficiently smooth, if
\begin{equation}\label{eq:cii}
C_{ij}(\delta_i)=
    \begin{cases}
\dfrac{4\kappa_i}{\pi \delta_i^4} &\quad   i=j,
\\[5mm]
\dfrac{\widetilde{C}_{ij}}{\delta_i^{4}} &\quad i \neq j,
\end{cases}
\end{equation}
then, in the local limit  $\delta_1\to0$ and $\delta_2\to0$ of vanishing horizons, the weak formulation of the nonlocal interface problem given in Proposition \ref{lem:weak-form} converges to its local counterpart \eqref{eq:weak_recovered} with constant $\kappa_1$ and $\kappa_2$. Moreover, the convergence rate is first order with respect to the horizon parameters.
\end{proposition}

Note that the constants and parameters in \eqref{eq:cii} associated with the one- and three-dimensional problems are different, but can be determined by following the same procedure presented for the two-dimensional case. In Section \ref{sec:num_res} we report examples of kernels for a one-dimensional problem.
We recall that for nonlocal volume-constrained problems such as the one in \eqref{eq:strong_nonloc} (i.e., in the absence of interfaces), the nonlocal operator converges to its local counterpart (the classical Laplacian) as $\mathcal O(\delta^2)$ and that the convergence of nonlocal solutions to the corresponding local ones is also of second order. We conjecture that also for the nonlocal interface problem, the rate of convergence of solutions is the same as that for the weak form, i.e., we conjecture that the convergence of nonlocal solutions to their corresponding local counterparts is also of first order. This conjecture is supported by the numerical results given in Section \ref{sec:num_res}.

Also note that, even though there is freedom in the choice of $\widetilde C_{ij}$, their values do affect the quality of the solution at the interface. In Section \ref{sec:num_res} we investigate the sensitivity of the nonlocal solution to different $\widetilde C_{ij}$.

\section{Numerical Results} \label{sec:num_res}

\begin{figure}[!t]
 \begin{center}
  \includegraphics[scale=0.9]{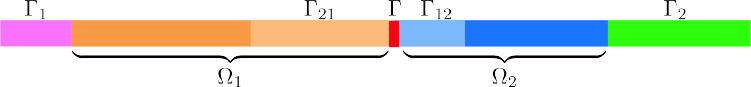}
   \caption{Domains configuration for the numerical results. Note that $\Gamma$ is a point in this one-dimensional setting.}
   \label{domains1D}
   \end{center}
\end{figure}

In this section we carry out a numerical investigation of the NLI theory developed in this paper. We first consider a one-dimensional setting, reported in Figure \ref{domains1D}, and then present preliminary tests in two dimensions, see  Figure \ref{2D}, top left, for the domain configuration. Note that, to simplify the numerical implementation, we consider an interface problem with a slightly different domain configuration compared to the one used in the previous sections; however, our theory applies to this configuration as well.

\subsection{One-dimensional problem}
In the section we illustrate the theoretical results presented in the previous sections and highlight some important features of the proposed approach.
\paragraph{Implementation details and problem setting}
The domains are discretized using an interface-fitted finite element grid $\mathcal{T}_h$ having $N_h$ nodes $x_i$ and $N_h-1$ elements of size $h$. We denote the node on the interface by $x_{i_\Gamma}$.
On each subdomain we approximate the nonlocal solution by a piecewise linear $C^0$ finite element space endowed with the standard Lagrangian nodal basis $\varphi_i$. To allow the nonlocal solution to develop a discontinuity on the interface we `` double-count'' the degree-of-freedom living on the interface node $x_{i_\Gamma}$.
Discretization of the nonlocal interface problem results in a $(N_h+1)\times (N_h+1)$ linear system of algebraic equations $A \mathbf{u} = \mathbf{f}$, where
\begin{equation}
\begin{aligned}\label{eq:linear_sys}
A_{ij} = \int_{\Omg \cup \widetilde\Gamma} \Big[\int_{\Omg \cup \widetilde\Gamma} \gamma(x, y) \Big(\varphi_j(x) - \varphi_j(y)\Big)\Big(\varphi_i(x) - \varphi_i(y)\Big)dy \Big] dx,
\qquad 
f_i = \int_\Omg f(x) \varphi_i(x) dx,
\end{aligned}
\end{equation}
for $i,j = 1, \ldots i_\Gamma, i_\Gamma+1,\ldots, N_h$. Integrals in the bilinear form are computed using a three-point Gauss quadrature (the same quadrature rule is used for error computation, descried later on).
In equation \eqref{eq:linear_sys}, $\Omg$ is an approximation of $\Omega_1 \cup \Omega_2$, and $\widetilde{\Gamma}$ of $\Gamma_1 \cup \Gamma_2$.
Let $\kappa_1$ and $\kappa_2$ be two positive constants describing the material properties in $\Omega_1$ and $\Omega_2$, respectively. By applying the theory of Section \ref{loc_lims} to a one-dimensional problem setting, we define the kernel in equation \eqref{eq:linear_sys} as follows.
\begin{equation}\label{eq:kernel1D}
\gamma(x,y) = \left\{
\begin{array}{ll}
\gamma_{11}(x,y) = \dfrac{3}{2} \, \dfrac{\kappa_1}{\delta_1^3} \mcX_{B_1(x)}(y)        & \;\; (x,y)\in \{\Omega_1\cup \Gamma_1 \times \Omega_1\cup \Gamma_1\} \\[6mm]
\gamma_{12}(x,y) = C_{12}(\delta_1)
\mcX_{B_1(x)}(y) & \;\; (x,y)\in \{\Omega_1\cup \Gamma_1 \times \Omega_2\cup \Gamma_2\}\\[6mm]
\gamma_{21}(x,y) =  C_{21}(\delta_2)
\mcX_{B_2(x)}(y) & \;\; (x,y)\in \{\Omega_2\cup \Gamma_2 \times \Omega_1\cup \Gamma_1\}  \\[6mm]
\gamma_{22}(x,y) = \dfrac{3}{2}\dfrac{\kappa_2}{\delta_2^3} \mcX_{B_2(x)}(y)        & \;\;(x,y)\in \{\Omega_2 \cup \Gamma_2 \times \Omega_2\cup \Gamma_2 \},
\end{array}\right.
\end{equation}
Four different kernels are considered, defined by the choice of $C_{12}$ and $C_{21}$.
\begin{align}
{\bf 1.} & \quad
C_{12}(\delta_1) = \dfrac{3}{2}\,\dfrac{\kappa_2}{\delta_1^3}, \quad
C_{21}(\delta_2) = \dfrac{3}{2} \, \dfrac{\kappa_1}{\delta_2^3}, 
\label{NewKernel}\\
{\bf 2.} & \quad
C_{12}(\delta_1) = \dfrac{3}{2}\,\dfrac{\kappa_1}{\delta_1^3}, \quad
C_{21}(\delta_2) = \dfrac{3}{2} \, \dfrac{\kappa_2}{\delta_2^3}.
\label{OldKernel1}\\
{\bf 3.} & \quad
C_{12}(\delta_1) = C_{21}(\delta_2)= \dfrac{3}{4} \, \Big(\dfrac{\kappa_1}{\delta_1^3}+\dfrac{\kappa_2}{\delta_2^3}\Big).
\label{OldKernel2}\\
{\bf 4.} & \quad
C_{12}(\delta_1) = \dfrac{3}{4} \, \Big(\dfrac{\kappa_1}{\delta_1^3}+\dfrac{\kappa_2}{\delta_1^3}\Big), \quad
C_{21}(\delta_2) = \dfrac{3}{4} \, \Big(\dfrac{\kappa_1}{\delta_2^3}+\dfrac{\kappa_2}{\delta_2^3}\Big).
\label{OldKernel2Modified}
\end{align}
The values of $C_{11}$ and $C_{22}$ are defined as in equation \eqref{eq:kernel1D}.
Unless otherwise stated, we let $\kappa_1=1$ and $\kappa_2=3$ in equation and consider the domains $\Gamma_1 = [-\delta_1 - 0.5, -0.5]$, $\Omega_1 = (-0.5,0)$, $\Omega_2 = (0,0.5)$ and $\Gamma_2=[0.5,0.5+\delta_2$], so that $\overline{\Omega} \cup \widetilde{\Gamma} = [-\delta_1-0.5,0.5+\delta_2]$. The forcing term is constant and such that $f(x)=1$ in $\Omega$. The volume constraints for the nonlocal problem are
\begin{align}\label{vol_constr_1D}
g_1(x) = \dfrac{1}{16}-\dfrac{1}{8}\,x-\dfrac{1}{2}\,x^2, \qquad g_2(x) = \dfrac{1}{16}-\dfrac{1}{24}\,x-\dfrac{1}{6}\,x^2.
\end{align}
In our study of the convergence to the local limits, also referred to as $\delta$-convergence, we consider the following local interface problem
\begin{align}\label{poissonSys}
\begin{cases}
- \kappa_1 u''_1(x) = f_1(x) \quad x\in(-0.5,0) \\[2mm]
- \kappa_2 u''_2(x) = f_2(x) \quad x\in(0,0.5) \\[2mm]
u_1(-0.5)=0, \quad u_2(0.5)=0, \quad u_1(0)=u_2(0)\\[2mm]
\kappa_1 u'_1(0) = \kappa_2 u'_2(0),
\end{cases} 
\end{align}
where, consistently with the nonlocal problem, $\kappa_1=1$, $\kappa_2=3$ and $f_1=f_2=1$. The analytic solution $u_L$ of the above problem is
\begin{align}\label{u_L}
u_L(x) = 
\begin{cases}
u_1(x) = g_1(x) & x \in (-0.5,0) \\
u_2(x) = g_2(x) & x \in (0,0.5),
\end{cases}
\end{align}
where $g_1$ and $g_2$ are the volume constraints of the nonlocal problem, defined in equation \eqref{vol_constr_1D}.

\paragraph{Convergence to the local limit}
To assess $\delta$-convergence, we consider a fine mesh with fixed size $h=2^{-12}$ and progressively halve $\delta_1$ and $\delta_2$, as $\|u_{N,h} - u_{L} \|_{L^2(\Omega\cup\widetilde\Gamma)}$ is monitored. Here, $u_{N,h}$ is the finite element nonlocal solution associated with the grid of size $h$. Since, in correspondence of the fine grid, the error $\|u_{L,h} - u_{L} \|_{L^2(\Omega)}$ between the numerical solution $u_{L,h}$ of problem \eqref{poissonSys} and the analytic solution $u_L$ is of order $10^{-9}$, we use $u_L$ in place of $u_{L,h}$ because we expect the discretization error to be negligible compared to the $\delta$ error. Note that from now on we drop the dependence of the norms on the domain. 

\begin{figure}[!t]
 \begin{center}
  \includegraphics[scale=0.45]{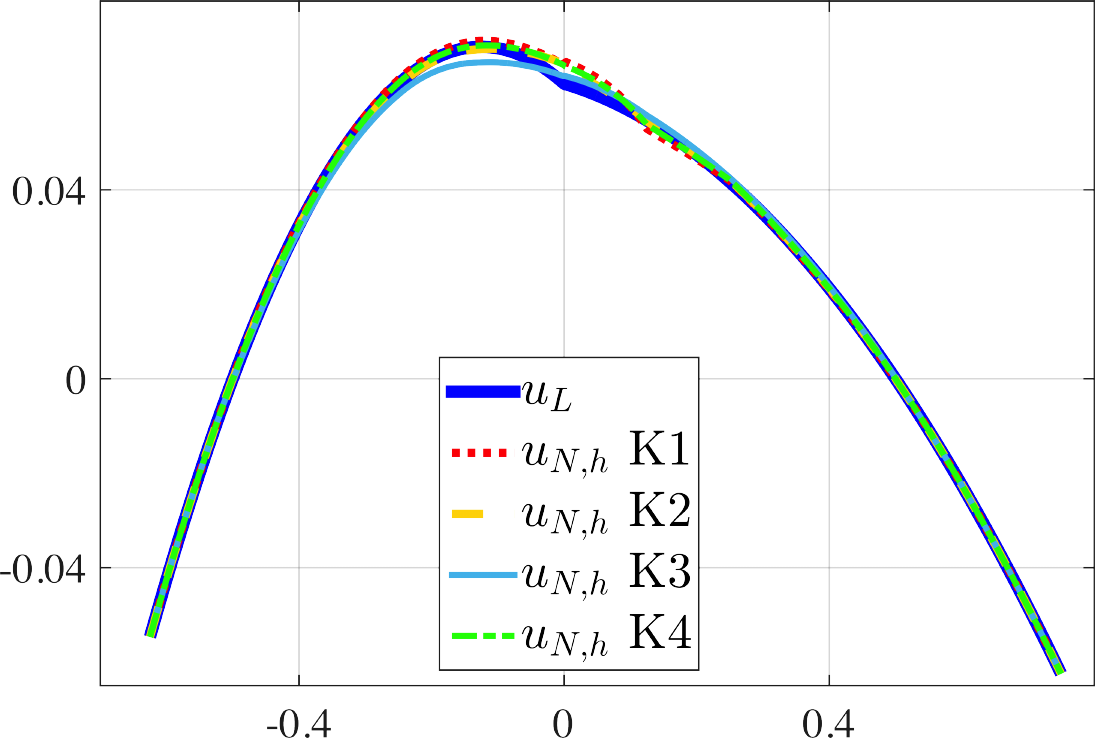}
    \includegraphics[scale=0.45]{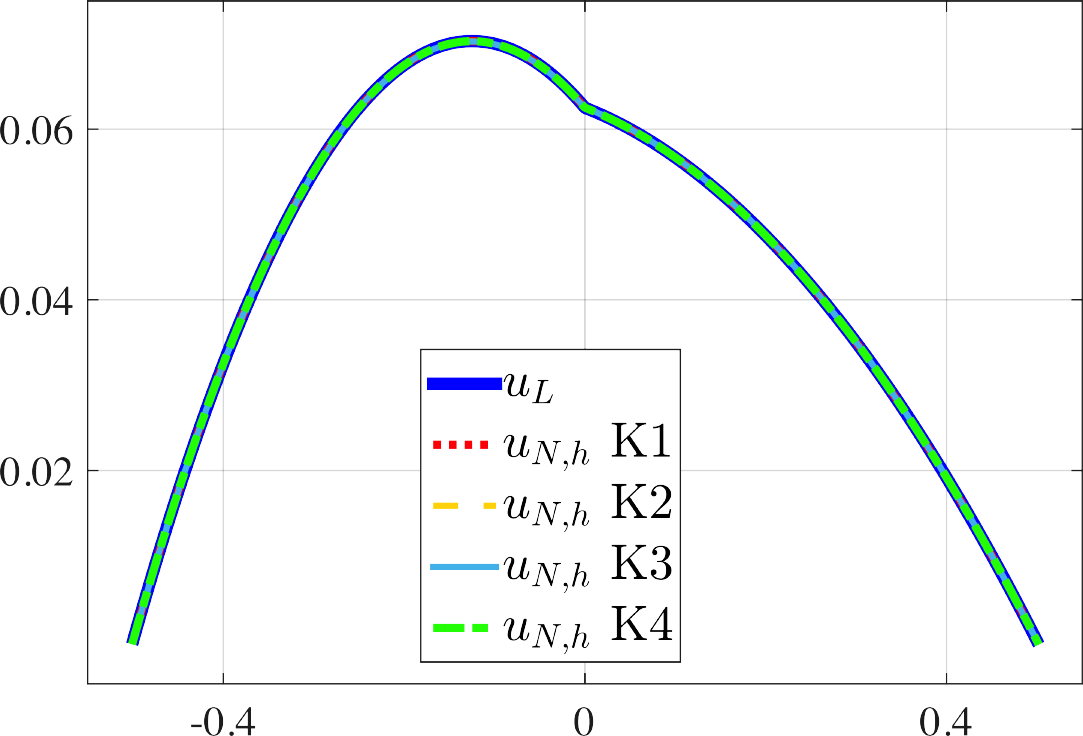}\\[2mm]
    \includegraphics[scale=0.56]{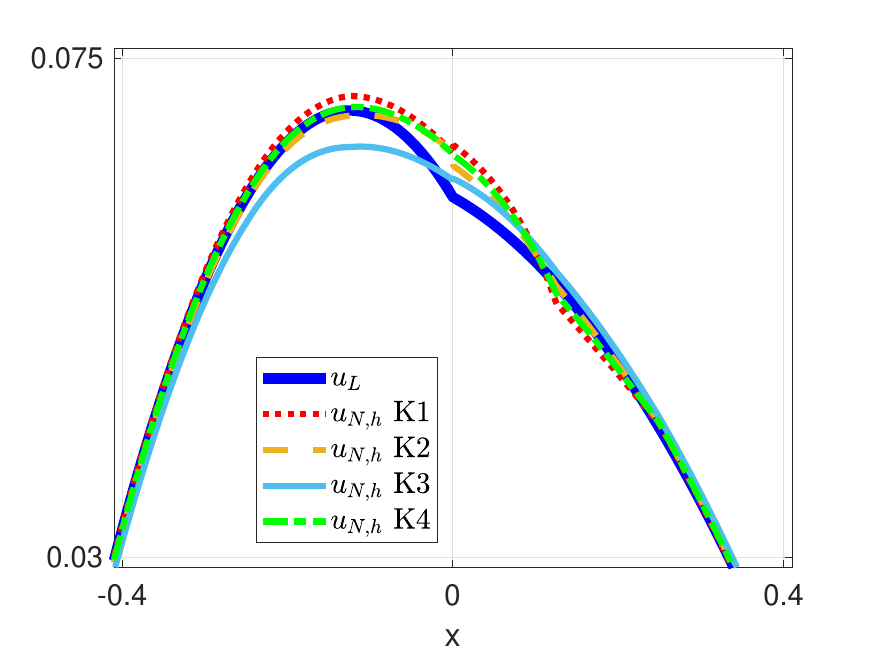}
    \includegraphics[scale=0.56]{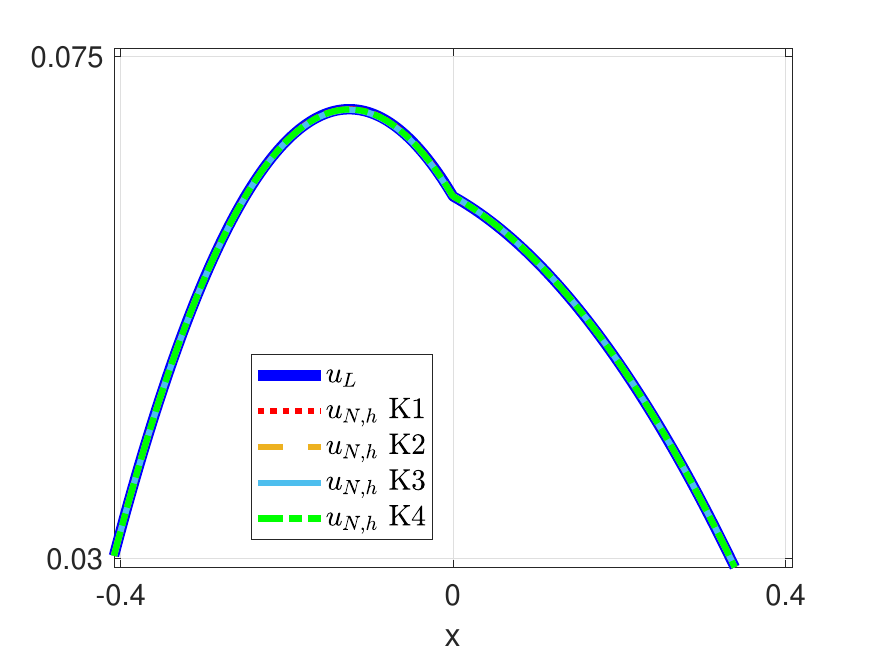}
   \caption{Left: comparison of nonlocal solution and exact local for different kernels with $\delta_1=2^{-10}, \delta_2=2^{-9}, h=2^{-12}.$ (zoom at the bottom). Right: comparison of nonlocal solution and exact local for different kernels with $\delta_1=2^{-3}, \delta_2=2^{-2}, h=2^{-12}.$ (zoom at the bottom).}
   \label{allkernels}
   \end{center}
\end{figure}

Results are reported in Table \ref{Tab1}. We see that the error between the nonlocal and the local solution goes to zero with a first order convergence for all kernels, showing how the local model can be recovered from the nonlocal model. Note that this confirms the conjecture made in Section \ref{loc_lims} that solutions of the nonlocal interface problems converge to its local counterpart as $\mathcal O(\delta_i)$.
In Figure \ref{allkernels}, we show the nonlocal solution obtained with the different kernels compared to the local exact $u_L$. On the left, we display solutions associated with relatively large values of $\delta_1$ and $\delta_2$, namely $\delta_1=2^{-3}$ and $\delta_2=2^{-2}$. 

As anticipated in Section \ref{loc_lims}, different $\widetilde C_{ij}$ yield a different behavior of the nonlocal solution across the interface. As an example, we note that the nonlocal solution obtained using Kernel 3 is almost insensitive to the presence of the interface. This behavior is likely due to the fact that $C_{12}=C_{21}$, i.e. the symmetric nature of the kernel prevents capturing the discontinuities in model parameters.

In Figure \ref{allkernels} (right), we show nonlocal solutions associated with relatively small values of $\delta_1$ and $\delta_2$, namely $\delta_1=2^{-10}$ and $\delta_2=2^{-9}$. This figure is meant to give visual proof that all kernels provide a nonlocal solution that converges to the local exact as $\delta_1$ and $\delta_2$ approach zero.
\begin {table}[ht]
\setlength\tabcolsep{3.25pt} 
\begin{center} 
\small
\begin{tabular}{cc|cc|cc|cc|cc}
& & \hspace{1cm}Kernel 1&  & \hspace{1cm}Kernel 2 & & \hspace{1cm}Kernel 3 & & \hspace{1cm}Kernel 4 & \\ \hline
$\delta_1$ & $\delta_2$  & $\|u_{N,h} - u_{L} \|_{L^2}$ & order & $\|u_{N,h} - u_{L} \|_{L^2}$ & order & $\|u_{N,h} - u_{L} \|_{L^2}$ & order & $\|u_{N,h} - u_{L} \|_{L^2}$ & order \\ \hline
$ 2^{-5} $& $  2^{-4}$  & 1.62e$-04$   &  --    &  3.86e$-04$  &  --   & 7.72e$-04$ & --    & 2.61e$-04$  & --\\  
$ 2^{-6} $& $  2^{-5}$  & 6.69e$-05$   &  1.28  &  2.19e$-04$  &  0.82 & 4.22e$-04$ & 0.87  & 1.45e$-04$  & 0.84 \\  
$ 2^{-7} $& $  2^{-6}$  & 3.11e$-05$   &  1.11  &  1.16e$-04$  &  0.91 & 2.20e$-04$ & 0.94  & 7.72e$-05$   & 0.91  \\  
$ 2^{-8} $& $  2^{-7}$  & 1.52e$-05$   &  1.04  &  6.01e$-05$  &  0.95 & 1.12e$-04$ & 0.97  & 3.98e$-05$   & 0.95  \\  
$ 2^{-9} $& $  2^{-8}$  & 7.52e$-06$   &  1.01  &  3.05e$-05$  &  0.98 & 5.68e$-05$ & 0.98  & 2.02e$-05$  & 0.98  \\  
$ 2^{-10}$& $ 2^{-9}$ & 3.75e$-06$   &  1.00  &  1.54e$-05$  &  0.99 & 2.86e$-05$ & 0.99  & 1.02e$-05$  & 0.99 \\
\end{tabular}
\end{center}
\caption{{\it One-dimensional problem}. Errors with respect to the local solution for decreasing values of $\delta_1$ and $\delta_2$.}
\label{Tab1}
\end{table}
\normalsize

\begin {table}[ht]
\setlength\tabcolsep{3.25pt} 
\begin{center} 
\small
\begin{tabular}{c|cc|cc|cc|cc} 
& \hspace{1cm}Kernel 1 & & \hspace{1cm} Kernel 2 & & \hspace{1cm}  Kernel 3 & & \hspace{1cm} Kernel 4 & \\ \hline
$h$      & $\|u_{N,h} - u_{N,h_f} \|_{L^2}$ & order   & $\|u_{N,h} - u_{N,h_f} \|_{L^2}$ & order 
 & $\|u_{N,h} - u_{N,h_f} \|_{L^2}$ & order 
  & $\|u_{N,h} - u_{N,h_f} \|_{L^2}$ & order \\ \hline 
    $2^{-5}$ & 6.58e$-05$ & --   & 5.86e$-05$& --   &  5.79e$-05$& --   & 6.07e$-05$  &--    \\
    $2^{-6}$ & 1.63e$-05$ & 2.01 & 1.36e$-05$& 2.10 &  1.32e$-05$& 2.13 & 1.44e$-05$ & 2.07\\
    $2^{-7}$ & 3.94e$-06$ & 2.05 & 3.33e$-06$& 2.03 &  4.08e$-06$& 1.69 & 3.43e$-06$ &2.07\\ 
    $2^{-8}$ & 9.49e$-07$ & 2.05 & 1.18e$-06$& 1.45 &  2.21e$-06$& 0.88 & 9.39e$-07$ &1.87\\ 
    $2^{-9}$ & 2.33e$-07$ & 2.02 & 6.77e$-07$& 0.80 &  1.40e$-06$& 0.65 & 4.25e$-07$ &1.14\\ 
\end{tabular}
\end{center}
\caption{{\it One-dimensional problem}. Errors with respect to a reference nonlocal solution for decreasing values of $h$.}
\label{Tab2}
\end{table}
\normalsize

\paragraph{Convergence with respect to the mesh size}
To investigate the $h$-convergence of the proposed finite element approximation, we let $u_{N,h_f}$ be a finite element nonlocal solution associated with a fine grid of size $h_f$, with $h_f \!\ll\! h$. 
For fixed values of $\delta_1$ and $\delta_2$, the $h$-convergence is assessed by progressively halving $h$ and monitoring $\|u_{N,h} - u_{N,h_f} \|_{L^2}$. 

Results are shown in Table \ref{Tab2}, for $\delta_1=2^{-5}, \delta_2=2^{-4}$, and $h_f=2^{-12}$.  Recall also that a double node is present at the interface, to allow a discontinuous nonlocal solution. From Table \ref{Tab2}, we see that the expected quadratic order of convergence given by the use of linear finite elements is obtained only for Kernel 1, whereas the other kernels display a rapidly deteriorating rate. In light of Table \ref{Tab2}, from now on we only consider Kernel 1, given in equation \eqref{NewKernel}\footnote{We are not able to theoretically explain the deterioration of convergence of Kernels 2--4; we believe a rigorous analysis of the numerical finite element error is needed; this is part of our current work.}.

\begin{figure}[!t]
\begin{center}
\begin{tabular}{cc}
\includegraphics[height=2in]{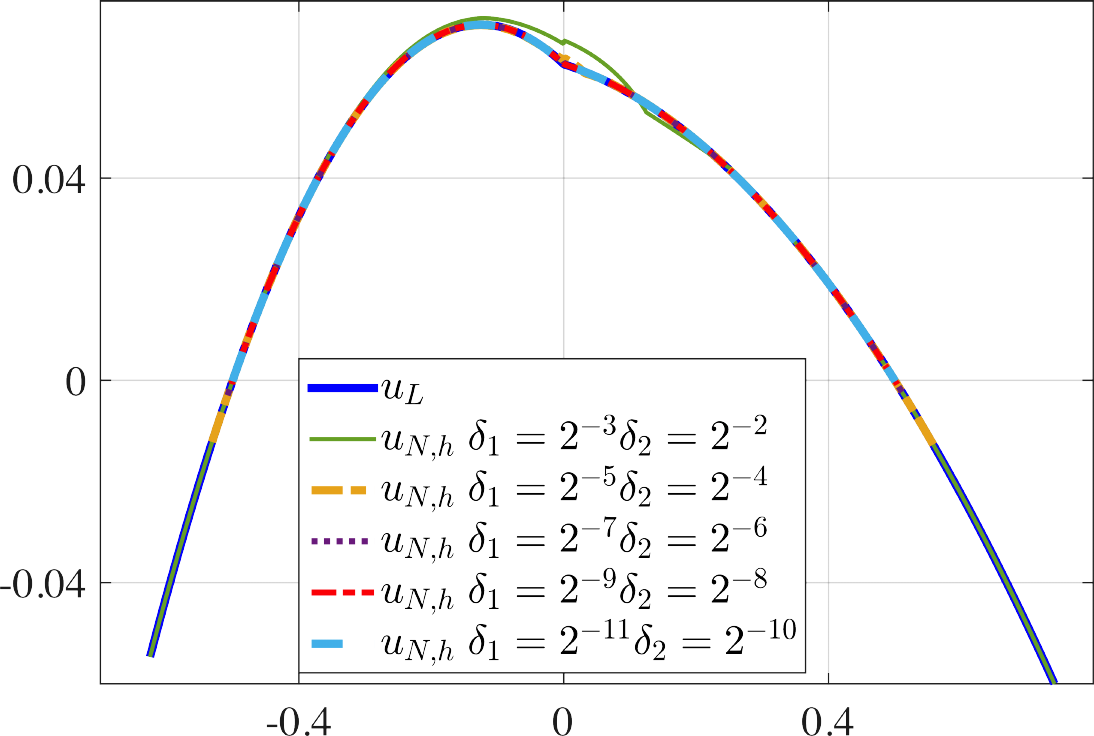} &
\hspace{.4cm}
\includegraphics[scale=0.42]{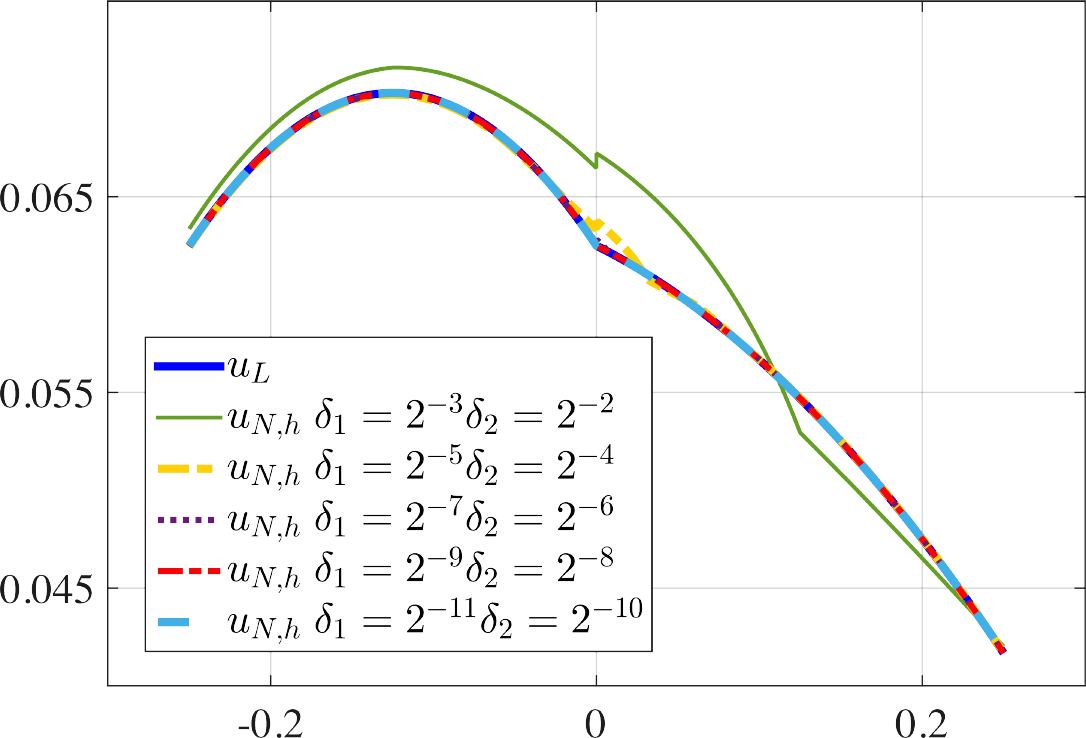}\\
\includegraphics[scale=0.42]{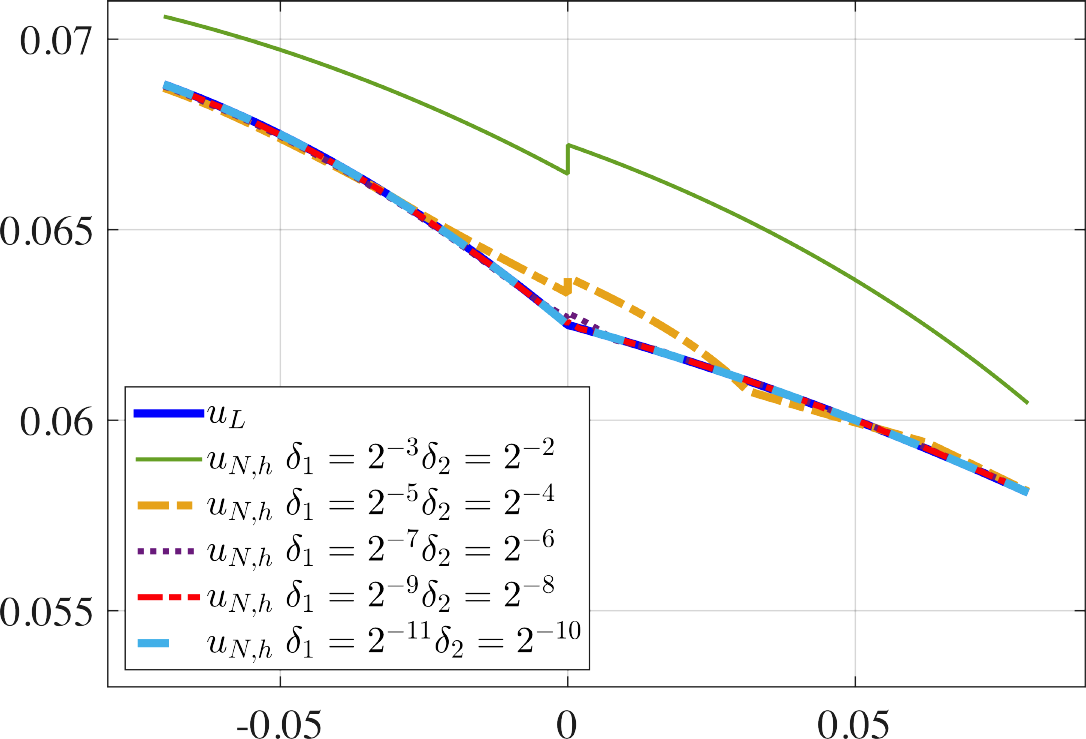} &
\hspace{.4cm}
\includegraphics[scale=0.42]{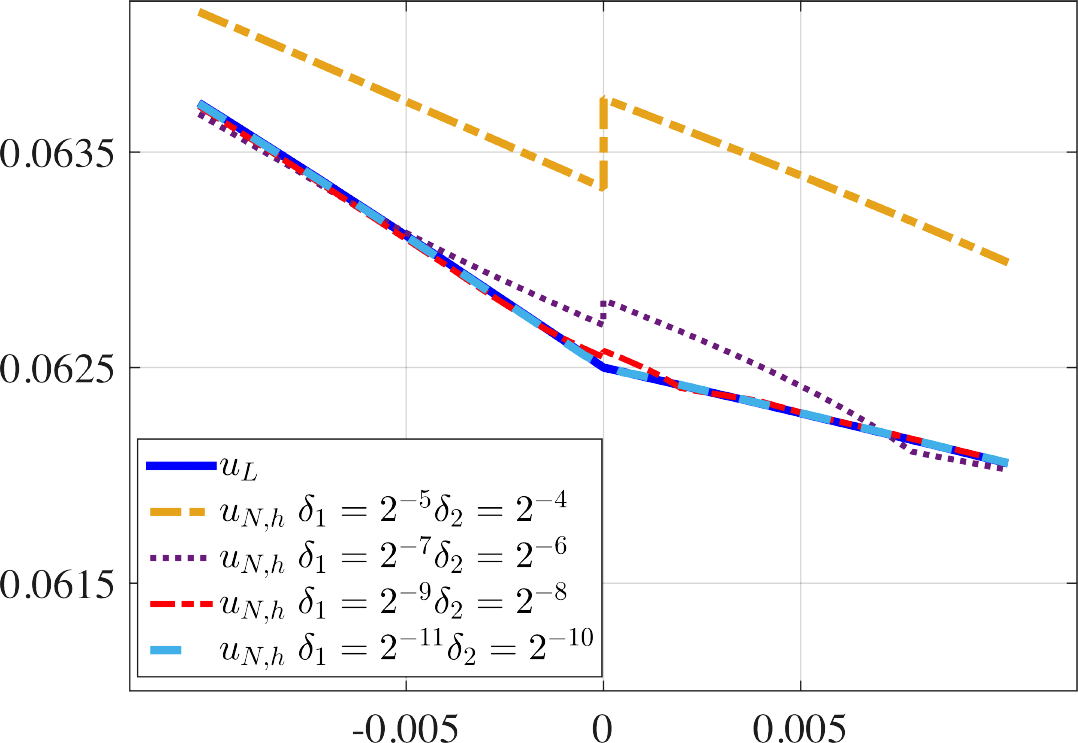}
\end{tabular}
\caption{Numerical nonlocal solutions compared to the local exact solution in equation \eqref{u_L}. The pictures are progressively zoomed on the interface moving from left to right and from top to bottom.}
\label{pics}
\end{center}
\end{figure}

\paragraph{Behavior of the solution at the interface}
First, we consider the behavior of the solution at the interface as the nonlocal interactions vanish. In Figure \ref{pics}, a plot of the numerical nonlocal solution for different values of $\delta_1$ and $\delta_2$ is compared to the local solution in equation \eqref{u_L}. It can be see from the pictures that the nonlocal solution has a jump discontinuity at the interface, and that the magnitude of the jump approaches zero as $\delta_1$ and $\delta_2$ approach zero, as confirmed by the results in Table \ref{Tab1}.

In Table \ref{Tab3}, left, we investigate the behavior of the magnitude of the jump, i.e. the difference of the solution values at the double node that has been placed at the interface. 
For $h=2^{-12}$, and progressively smaller values of $\delta_1$ and $\delta_2$, we report the magnitude of the discontinuity; we observe that it approaches zero with first order convergence. This shows that the nonlocal solution starts as discontinuous, and as it converges to the local solution becomes continuous.

In Table \ref{Tab3}, right, we consider $\delta_1=2^{-5}$, $\delta_2=2^{-4}$ as $h$ is decreased. As expected, the magnitude of the discontinuity reaches a saturation value. This behavior is due to the fact that the discontinuity is intrinsically related to nonlocality and its magnitude depends on the values of $\delta_1$ and $\delta_2$.
\begin{table}[t]
\setlength\tabcolsep{3.25pt} 
\begin{center}
\small
\begin{tabular}{ll|cc}	

$\delta_1$ & $\delta_2$ & magnitude  & order  \\ \hline 
$2^{-5}$   & $2^{-4}$   & 4.15e$-04$ & --     \\ 
$2^{-6}$   & $2^{-5}$   & 2.25e$-04$ & 0.88   \\ 
$2^{-7}$   & $2^{-6}$   & 1.17e$-04$ & 0.94   \\ 
$2^{-8}$   & $2^{-7}$   & 5.95e$-05$ & 0.97   \\ 
$2^{-9}$   & $2^{-8}$   & 3.00e$-05$ & 0.99   \\ 
$2^{-10}$  & $2^{-9}$   & 1.51e$-05$ & 0.99   \\ 
\end{tabular}
\hspace{1cm}
\begin{tabular}{l|cc}
$h$       & magnitude  & order       \\ \hline  
$2^{-5}$  & 6.50e$-04$ & --          \\ 
$2^{-6}$  & 4.23e$-04$ & 6.20e$-01$  \\ 
$2^{-7}$  & 4.17e$-04$ & 1.99e$-02$  \\ 
$2^{-8}$  & 4.15e$-04$ & 6.20e$-03$  \\ 
$2^{-9}$  & 4.15e$-04$ & 1.00e$-04$  \\ 
$2^{-10}$ & 4.15e$-04$ & 3.00e$-04$  \\ 
$2^{-11}$ & 4.15e$-04$ & 0.00e$+00$  \\ 
\end{tabular}
\end{center}
\caption{{\it One-dimensional problem}. Left: magnitude of the jump discontinuity at the interface of the nonlocal solution, for fixed $h$ and decreasing $\delta_1$ and $\delta_2$. Right: magnitude of the jump discontinuity at the interface of the nonlocal solution, for fixed $\delta_1$ and $\delta_2$ and decreasing $h$.}
\label{Tab3}
\end{table}
\normalsize

\begin{figure}[!t]
\begin{center}
\includegraphics[scale=0.3]{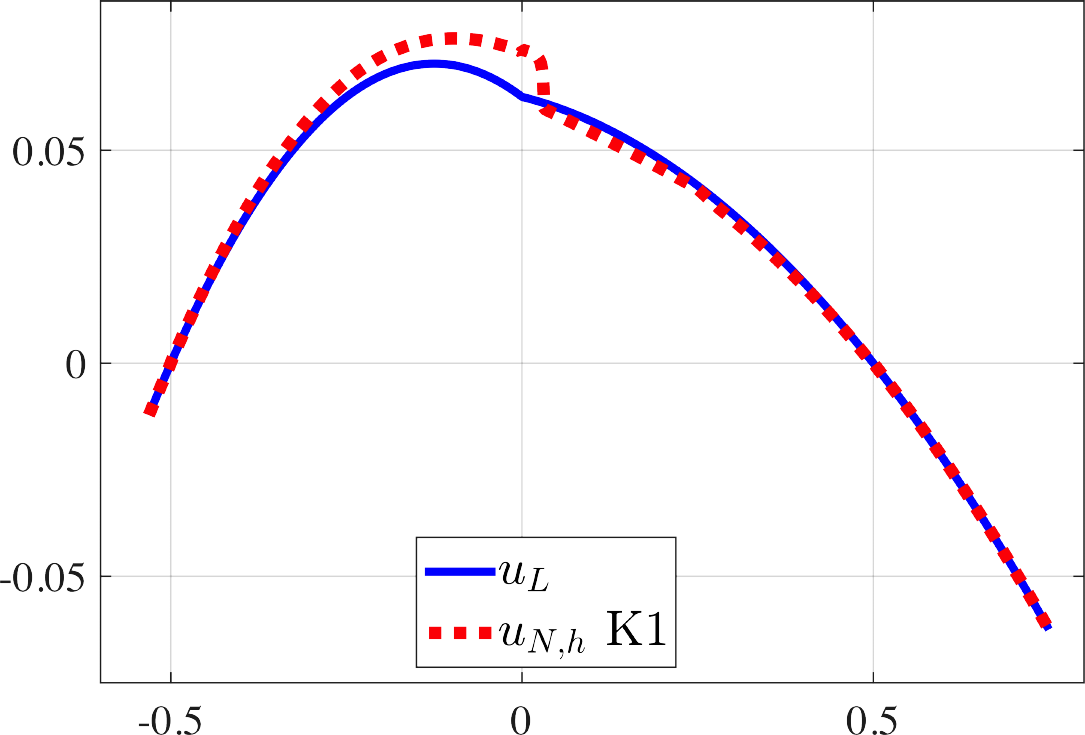}
\includegraphics[scale=0.3]{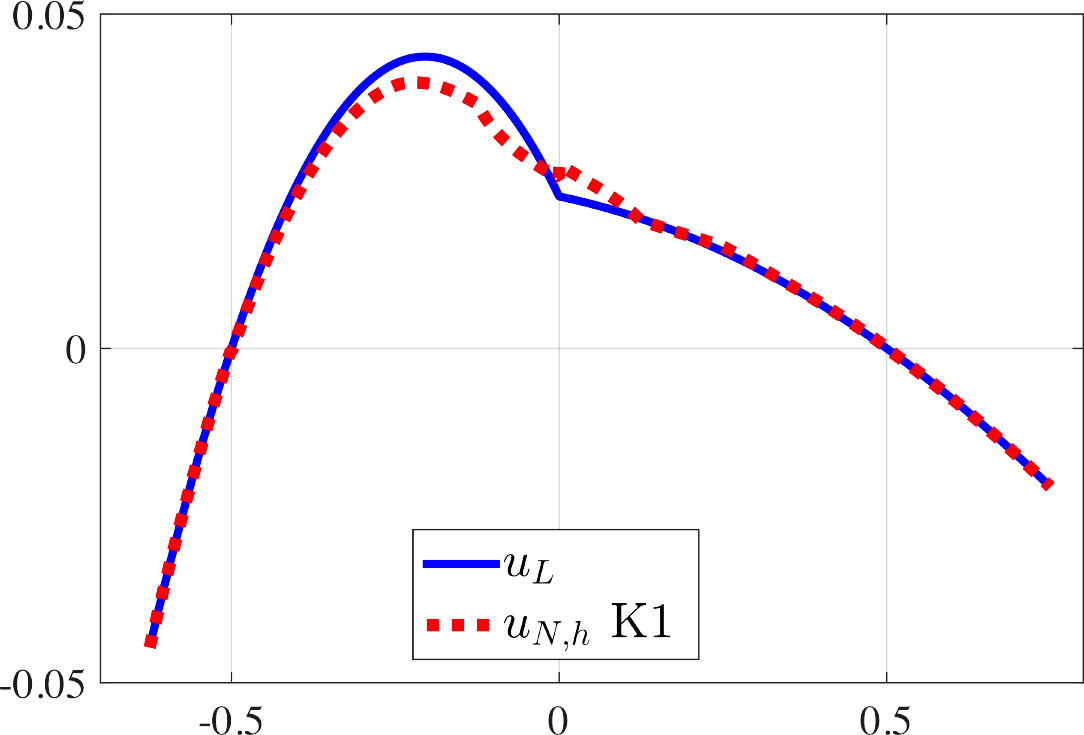}
\includegraphics[scale=0.3]{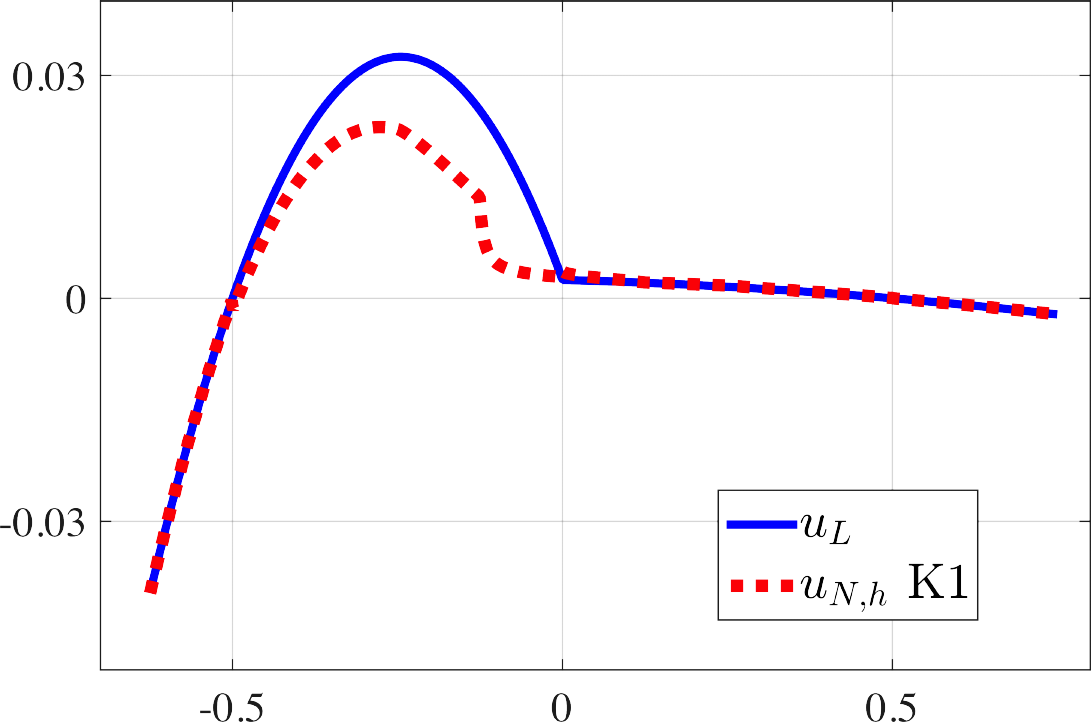}
\caption{Nonlocal solution for $h=2^{-12}$ and left: $\kappa_1=1$, $\kappa_2=3$, $\delta_1=2^{-5}$, $\delta_2=2^{-2}$; center: $\kappa_1=1$, $\kappa_2=10$, $\delta_1=2^{-3}$, $\delta_2=2^{-2}$; right: $\kappa_1=1$, $\kappa_2=100$, $\delta_1=2^{-3}$, $\delta_2=2^{-2}$.}
\label{mix}
\end{center}
\end{figure}

\paragraph{Sensitivity to the parameters $\kappa_i$}
To complete the one-dimensional investigation, we also report pictures of nonlocal solutions obtained with a larger difference between the values of $\kappa_1$  and $\kappa_2$, or between the values of $\delta_1$ and $\delta_2$, see Figure \ref{mix}. Note that, while the values of $\kappa_i$ change the solution profile, the convergence behavior with respect to $\delta_i$ and $h$ is not affected, as it is independent of them.

\subsection{Two-dimensional problem}
\begin{figure}[!t]
\begin{center}
\begin{tabular}{cc}
\includegraphics[height=2in]{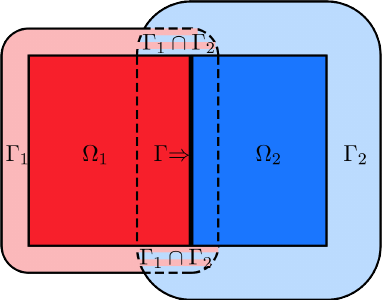} & \hspace{.5cm}
\includegraphics[scale=0.5]{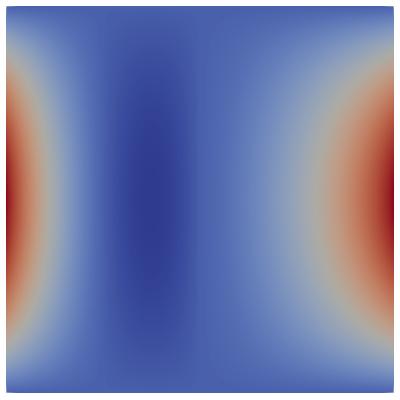} \\[0.3cm]
\includegraphics[scale=0.46]{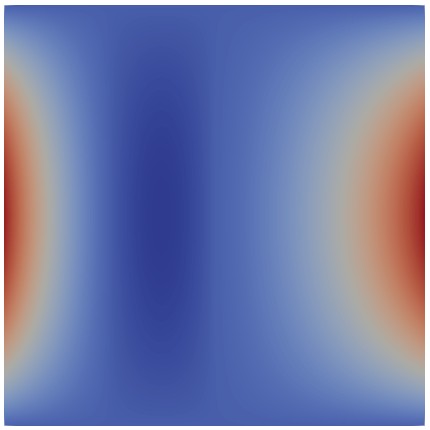} & \hspace{.5cm}
\includegraphics[scale=0.5]{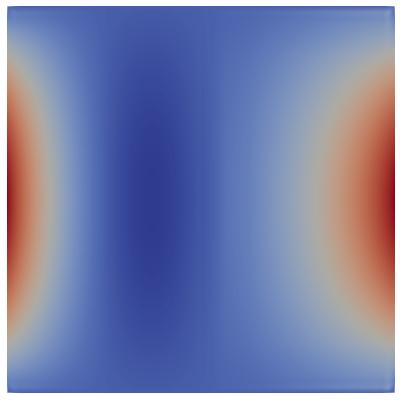} 
\end{tabular}
\caption{Top left: domains configuration for the 2D numerical tests. Top Right: numerical local solution in $\Omega_1 \cup \Omega_2$. Bottom left: numerical nonlocal solution in $\Omega_1 \cup \Omega_2$ with $(\delta_1,\delta_2) = (2^{-8},2^{-7})$. Bottom right: numerical nonlocal solution in $\Omega_1 \cup \Omega_2$ with $(\delta_1,\delta_2) = (2^{-4},2^{-3})$.}
\label{2D}
\end{center}
\end{figure}

In this section we show the applicability of our strategy to higher-dimensional problems and illustrate the theoretical results in Section \ref{loc_lims}. We refer to the configuration in Figure \ref{2D}, top left; also in this case, the domains are discretized using an interface-fitted finite element grid of size $h$. On each subdomain the nonlocal solution is a piecewise linear finite element approximation. Double edges and nodes are placed on the interface to allow for a discontinuous nonlocal solution across the interface. 

\paragraph{Problem setting}
We consider $\Gamma_1 \cup \Omega_1 = [-\delta_1-0.5, \delta_1]^2$,  and $\Omega_2 \cup \Gamma_2 = [-\delta_2,0.5+\delta_2]^2$. Let $\xb=(x_1,x_2)$, the nonlocal volume constraints on the interaction domains are defined as

\begin{align}
g_1(\xb) &= 
\begin{cases}
1/16  & \quad\;\; |x_2| > 0.5\\
1/16+ (-1/8\,x_1-1/2 \,x_1^2)\Big(x_2^2 - 1/4\Big) & \quad\;\; \mbox{otherwise}
\end{cases}\\[2mm]
g_2(\xb) &= 
\begin{cases}
1/16  & \,\, |x_2| > 0.5\\
1/16 + \Big(-1/24\,x_1-1/6\,x_1^2\Big)\Big(x_2^2 - 1/4\Big) & \,\, \mbox{otherwise},
\end{cases}
\end{align}
and, for $\kappa_1=1$ and $\kappa_2=3$, the forcing term is defined as $f(\xb)= - \kappa_i \Delta g_i$ for $\xb \in \Omega_i$.
In the analysis of the $\delta$-convergence, we consider the following local problem 
\begin{align}\label{poissonSys2D_2}
\begin{cases}
-\kappa_1 \Delta u_1 = f  & \xb\in\Omega_1 \\[1mm]
-\kappa_2 \Delta u_2 = f  & \xb\in\Omega_2 \\[1mm]
u_1 = g_1                 & \xb\in\partial \Omega_1\setminus\Gamma\\[1mm]
u_2 = g_2                 & \xb\in\partial \Omega_2\setminus\Gamma\\[1mm]
u_1(\xb)=u_2(\xb)         & \xb\in\Gamma\\[1mm]
\kappa_1\dfrac{\partial{u_1}}{\partial x_1} = \kappa_2\dfrac{\partial{u_2}}{\partial x_1} & \xb\in\Gamma,
\end{cases} 
\end{align}
whose analytic solution $u_L$, reported in Figure \ref{2D}, top right, is given by $u_L(\bm{x}) = g_i(\xb)$ for $\xb\in \Omega_i$. The nonlocal kernel $\gamma$ is chosen in accordance to the results in Section \ref{loc_lims}.

\paragraph{Convergence to the local limit}
We conduct the same analysis of the previous section and analyze the convergence of the finite element nonlocal solution to the approximate solution of problem \eqref{poissonSys2D_2} as $\delta_1$ and $\delta_2$ approach zero. Note that in this case, we do not use the analytic local solution as the finite element grid is not fine enough to make the discretization error negligible. We denote by $u_{L,h}$ the local finite element solution on a grid of size $h$.

In Table \ref{Tab2D} we report values of $\|u_{N,h}-u_{L,h} \|_{L^2}$ for $h=2^{-8}$ as the interaction radii approach zero. Also in this case, we observe a first order convergence to the local solution. In Figure \ref{2D}, bottom, we report $u_{N,h}$ in $\Omega_1\cup\Omega_2$ for $(\delta_1,\delta_2)=(2^{-8},2^{-7})$, bottom left, and $(2^{-4},2^{-3})$, bottom right.
\begin {table}[h]
\setlength\tabcolsep{3.25pt} 
\begin{center}                                          \begin{tabular}{ll|cc} 	
$\delta_1$ & $\delta_2$ &  $\|u_{N,h}-u_{L,h} \|_{L^2}$ & order \\ \hline 
$2^{-4}$   & $2^{-3}$   &  3.06e$-04$  & --     \\ 
$2^{-5}$   & $2^{-4}$   &  1.42e$-04$  &  1.10  \\ 
$2^{-6}$   & $2^{-5}$   &  6.64e$-05$  &  1.10  \\ 
$2^{-7}$   & $2^{-6}$   &  2.94e$-05$  &  1.18  \\ 
$2^{-8}$   & $2^{-7}$   &  1.17e$-05$  &  1.32  \\ 
\end{tabular}
\end{center}
\caption{{\it Two-dimensional problem}. Difference between numerical nonlocal and local solutions as the interaction radii approach zero on a mesh of size $h=2^{-8}$.}
\label{Tab2D} 
\end{table}

\section{Conclusion and Perspectives}
\label{sec:conclusions}
We developed and demonstrated a new mathematical formulation for nonlocal interface problems (NLI), which extends the classical local interface theory to the non-local setting. 
Our theory, based on the minimization of the energy of the nonlocal system, provides a rigorous and physically consistent nonlocal counterpart of the classical theory, and fills a longstanding theoretical gap in the formulation of nonlocal interface problems. 
 An important feature of our approach is that as the extent of the nonlocal interactions vanishes, the solution of the nonlocal interface problem converges to the one of corresponding local problem; we refer to this property as {\it physical consistency.} Convergence properties of our formulation, both with respect to the nonlocal parameter and the discretization size, are illustrated by several one-dimensional experiments. Furthermore, a two-dimensional experiment shows the applicability of our strategy in higher dimensions and represents a promising preliminary result towards realistic simulations.

In this work we focused on nonlocal generalizations of perfect interface conditions. Subsequent work will address application of the NLI theory to the design of efficient nonlocal domain decomposition solvers and its extension to imperfect interfaces that occur in important applications such as fracture mechanics and problems with interfacial thermal conductance \cite{Mahan_09_PRB}. 
In particular, our objectives include  1) Consistent  NLI formulations for problems in which the solution has a prescribed jump at the interface; 2) extension of the numerical tests to more complex geometries; and 3) extension of NLI to singular kernels, that are characteristic of fracture problems and subsurface flow applications.

\noindent\section*{Acknowledgments}

\noindent This work was supported by the Sandia National Laboratories (SNL) Laboratory-directed Research and Development (LDRD) program, and the U.S. Department of Energy, Office of Science, Office of Advanced Scientific Computing Research under Award Number DE-SC-0000230927 and under the Collaboratory on Mathematics and Physics-Informed Learning Machines for Multiscale and Multiphysics Problems (PhILMs) project. SNL is a multimission laboratory managed and operated by National Technology and Engineering Solutions of Sandia, LLC., a wholly owned subsidiary of Honeywell International, Inc., for the U.S. Department of Energy's National Nuclear Security Administration under contract DE-NA-0003525. 

\smallskip
This paper describes objective technical results and analysis. Any subjective views or opinions that might be expressed in the paper do not necessarily represent the views of the U.S. Department of Energy or the United States Government. Report number SAND2020-4918.

\noindent \section*{In memoriam}
 
\noindent This paper is dedicated to the memory of Dr. Douglas Nelson Woods (January 11 1985 - September 11 2019),
promising young scientist and post-doctoral research fellow at Los Alamos National Laboratory.  
Our thoughts and wishes go to his wife Jessica, to his parents Susan and Tom, to his sister Rebecca and to his brother Chris, whom he left behind.


\appendix


\section{Proof of Proposition \ref{lem:weak-form}}\label{sec:weak-form-proof}
The Euler-Lagrange equation corresponding to the Minimization Principle \ref{eq:minprin-NLI} is given by
$$
\begin{aligned}
 &\int\limits_{\Omega_1 \cup\Gami}\int\limits_{\Omega_1\cup\Gami}
    (u_1(\xb)-u_1(\yb))(v_1(\xb)-v_1(\yb)) \gamma(\xb,\yb)d\yb\,d\xb
\\&\qquad+ \int\limits_{\Omega_1}\int\limits_{\Omg_2}
(u_1(\xb)-u_1(\yb))(v_1(\xb)-v_1(\yb)) \gamma(\xb,\yb)d\yb\,d\xb
- \int\limits_{\Omega_1} f v_1\, d\xb \\[2mm]
&\qquad+\int\limits_{\Omg_2}\int\limits_{\Omg_2} (u_2(\xb)-u_2(\yb))(v_2(\xb)-v_2(\yb)) \gamma(\xb,\yb)d\yb\,d\xb
\\&\qquad+\int\limits_{\Omg_2}\int\limits_{\Omega_1} (u_2(\xb)-u_2(\yb))(v_2(\xb)-v_2(\yb)) \gamma(\xb,\yb)d\yb\,d\xb
- \int\limits_{\Omg_2} f v_2\, d\xb = 0,
\end{aligned}
$$
for all $(v_1,v_2)\in W^c$ satisfying $v_1(\xb)=v_2(\xb) $ for $ \xb \in \Gamma^*$ and $v_1(\xb)=0$ for $ \xb \in\Gami$, where we have again used the assumption that points in $\Omg_2$ do not interact with points in $\Gami$. Substituting \eqref{eq:kernel}, \eqref{eq:cond_en_princ_nl}, and $v_1(\xb)=v_2(\xb) $ for $ \xb \in \Gamma^*$ results in
$$
\begin{aligned}
 &\int\limits_{\Omega_1 \cup\Gami}\int\limits_{\Omega_1\cup\Gami}
    (u_1(\xb)-u_1(\yb))(v_1(\xb)-v_1(\yb)) \gamma_{11}(\xb,\yb)d\yb\,d\xb
\\&\qquad+ \int\limits_{\Omega_1}\int\limits_{\Omg_2}
(u_1(\xb)-u_2(\yb))(v_1(\xb)-v_2(\yb)) \gamma_{12}(\xb,\yb)d\yb\,d\xb
- \int\limits_{\Omega_1} f v_1\, d\xb \\[2mm]
&\qquad+\int\limits_{\Omg_2}\int\limits_{\Omg_2} (u_2(\xb)-u_2(\yb))(v_2(\xb)-v_2(\yb)) \gamma_{22}(\xb,\yb)d\yb\,d\xb
\\&\qquad+\int\limits_{\Omg_2}\int\limits_{\Omega_1} (u_2(\xb)-u_1(\yb))(v_2(\xb)-v_1(\yb)) \gamma_{21}(\xb,\yb)d\yb\,d\xb
- \int\limits_{\Omg_2} f v_2\, d\xb = 0.
\end{aligned}
$$
Rearranging terms, we obtain
$$
\begin{aligned}
 &\int\limits_{\Omega_1 \cup\Gami}\int\limits_{\Omega_1\cup\Gami}
    (u_1(\xb)-u_1(\yb))(v_1(\xb)-v_1(\yb)) \gamma_{11}(\xb,\yb)d\yb\,d\xb
\\&\qquad+ \int\limits_{\Omega_1}\int\limits_{\Omg_2}
(u_1(\xb)-u_2(\yb))v_1(\xb) \gamma_{12}(\xb,\yb)d\yb\,d\xb
\\
&\qquad-\int\limits_{\Omg_2}\int\limits_{\Omega_1} (u_2(\xb)-u_1(\yb))v_1(\yb) \gamma_{21}(\xb,\yb)d\yb\,d\xb
- \int\limits_{\Omega_1} f v_1\, d\xb \\[2mm]
&+\int\limits_{\Omg_2}\int\limits_{\Omg_2} (u_2(\xb)-u_2(\yb))(v_2(\xb)-v_2(\yb)) \gamma_{22}(\xb,\yb)d\yb\,d\xb
\\&\qquad+\int\limits_{\Omg_2}\int\limits_{\Omega_1} (u_2(\xb)-u_1(\yb))v_2(\xb)) \gamma_{21}(\xb,\yb)d\yb\,d\xb
\\&\qquad- \int\limits_{\Omega_1}\int\limits_{\Omg_2}
(u_1(\xb)-u_2(\yb))v_2(\yb) \gamma_{12}(\xb,\yb)d\yb\,d\xb- \int\limits_{\Omg_2} f v_2\, d\xb = 0.
\end{aligned}
$$
Then, \eqref{eq:weak-omega1} and \eqref{eq:weak-omega2} follow because $v_1(\xb)$ for $\xb\in\Omg_1$ and $v_2(\xb)$ for $\xb\in\Omg_1$ can be independently chosen.


\section{Proof of Proposition \ref{lem:strong-form}}\label{sec:strong-form-proof}
We make use of the identity
\begin{align}\label{intparts}
\int_{D_1}\int_{D_2} \psi(\xb,\yb)\sigma(\yb)\gamma(\xb,\yb) d\yb d\xb  
= \int_{D_2}\sigma(\xb)\int_{D_1} \psi(\yb,\xb)\gamma(\yb,\xb) d\yb d\xb,
\end{align}
where $D_1$ and $D_2$ are two generic subsets of $\mathbb{R}^n$. 
For the first term in \eqref{eq:weak-omega1} we have, using \eqref{intparts}, $v_1(\xb)=0$ on $\Gami$, the first choice in \eqref{eq:kernel}, and the symmetry of $\gamma_{11}(\xb,\yb)$,
\begin{equation}\label{eq:weakA}
\begin{aligned}
& \int\limits_{\Omega_1 \cup\Gami}\int\limits_{\Omega_1\cup\Gami}
    (u_1(\xb)-u_1(\yb))(v_1(\xb)-v_1(\yb)) \gamma_{11}(\xb,\yb)d\yb\,d\xb
     \\[2mm]
  &= \int\limits_{\Omega_1\cup\Gami}v_1(\xb)\int\limits_{\Omega_1\cup\Gami}
    (u_1(\xb)-u_1(\yb))\gamma_{11}(\xb,\yb)\,d\yb\,d\xb
   -\int\limits_{\Omega_1\cup\Gami}v_1(\xb)\int\limits_{\Omega_1\cup\Gami}
    (u_1(\yb)-u_1(\xb))\gamma_{11}(\yb,\xb)\,d\yb\,d\xb \\[2mm]  
  &=-2\int\limits_{\Omega_1}v_1(\xb)\int\limits_{\Omega_1\cup\Gami}
    (u_1(\yb)-u_1(\xb))\gam_{11}(\xb,\yb)\,d\yb\,d\xb.
\end{aligned}
\end{equation}
For the third term in \eqref{eq:weak-omega1} we have,
\begin{equation}\label{eq:weakB}
\begin{aligned}
-\int\limits_{\Omg_2}\int\limits_{\Omega_1} (u_2(\xb)-u_1(\yb))v_1(\yb) \gamma_{21}(\xb,\yb)d\yb\,d\xb
=
-\int\limits_{\Omg_1}v_1(\xb)\int\limits_{\Omega_2} (u_2(\yb)-u_1(\xb)) \gamma_{21}(\xb,\yb)d\yb\,d\xb.
\end{aligned}
\end{equation}
Substituting \eqref{eq:weakA} and \eqref{eq:weakB} in \eqref{eq:weak-omega1} results in
$$
\begin{aligned}
&-2\int\limits_{\Omega_1}v_1(\xb)\int\limits_{\Omega_1\cup\Gami}
    (u_1(\yb)-u_1(\xb))\gam_{11}(\xb,\yb)\,d\yb\,d\xb
- \int\limits_{\Omega_1}v_1(\xb)\int\limits_{\Omg_2}
(u_2(\yb)-u_1(\xb)) \gamma_{12}(\xb,\yb)d\yb\,d\xb
\\&\qquad-\int\limits_{\Omg_1}v_1(\xb)\int\limits_{\Omega_2} (u_2(\yb)-u_1(\xb)) \gamma_{21}(\xb,\yb)d\yb\,d\xb
= \int\limits_{\Omega_1} f(\xb) v_1(\xb)\, d\xb.  
\end{aligned}
$$
Because $v_1(\xb )$ for $\xb\in\Omg_1$ is arbitrary, \eqref{eq:strong1}  follows. In a similar manner, \eqref{eq:strong2} is derived from \eqref{eq:weak-omega2}.

\section{Proof of Proposition \ref{lem:local-limit}}\label{sec:local-limit-proof}
Let us define
\begin{equation}
\begin{aligned}
    &NL(u_1,u_2,v_1,v_2,\gamma):= \int\limits_{\Omega_1 \cup\Gamma_1}\int\limits_{\Omega_1\cup\Gamma_1}
    (u_1(\xb)-u_1(\yb))(v_1(\xb)-v_1(\yb)) \gamma_{11}(\xb,\yb)d\yb\,d\xb\\
    &+ \int\limits_{\Omega_1\cup\Gamma_1}\int\limits_{\Omg_2}
(u_1(\xb)-u_1(\yb))(v_1(\xb)-v_1(\yb)) \gamma_{12}(\xb,\yb)d\yb\,d\xb 
+\int\limits_{\Omg_2}\int\limits_{\Omg_2} (u_2(\xb)-u_2(\yb))(v_2(\xb)-v_2(\yb)) \gamma_{22}(\xb,\yb)d\yb\,d\xb\\
&+\int\limits_{\Omg_2}\int\limits_{\Omega_1\cup\Gamma_1} (u_2(\xb)-u_2(\yb))(v_2(\xb)-v_2(\yb)) \gamma_{21}(\xb,\yb)d\yb\,d\xb - \int\limits_{\Omega_1} f_1 v_1\, d\xb - \int\limits_{\Omg_2} f_2 v_2. d\xb. 
\end{aligned}
\end{equation}
\begin{equation}
    \begin{aligned}
  L(u_1,u_2,v_1,v_2):=  \int_{\Omega_1} \kappa_1\nabla u_1(\xb)\cdot \nabla v_1(\xb)d\xb
+\int_{\Omega_2} \kappa_2\nabla u_2(\xb)\cdot \nabla v_2(\xb)d\xb  
-\int_{\Omega_1} f_1v_1d\xb 
-\int_{\Omega_2} f_2v_2d\xb.
    \end{aligned}
\end{equation}
Hence, we have
\begin{align}\label{eq:terms}
& NL(u_1,u_2,v_1,v_2,\gamma) - L(u_1,u_2,v_1,v_2)  \\[5mm]
&= \int\limits_{\Omega_1 \cup\Gamma_1}\int\limits_{\Omega_1\cup\Gamma_1}
    (u_1(\xb)-u_1(\yb))(v_1(\xb)-v_1(\yb)) \gamma_{11}(\xb,\yb)d\yb\,d\xb
+ \int\limits_{\Omega_1\cup\Gamma_1}\int\limits_{\Omg_2}
(u_1(\xb)-u_1(\yb))(v_1(\xb)-v_1(\yb)) \gamma_{12}(\xb,\yb)d\yb\,d\xb \nonumber \\
&+\int\limits_{\Omg_2}\int\limits_{\Omg_2} (u_2(\xb)-u_2(\yb))(v_2(\xb)-v_2(\yb)) \gamma_{22}(\xb,\yb)d\yb\,d\xb
+\int\limits_{\Omg_2}\int\limits_{\Omega_1\cup\Gamma_1} (u_2(\xb)-u_2(\yb))(v_2(\xb)-v_2(\yb)) \gamma_{21}(\xb,\yb)d\yb\,d\xb \nonumber\\
&- \int_{\Omega_1} \kappa_1\nabla u_1(\xb)\cdot \nabla v_1(\xb)d\xb  
-\int_{\Omega_2} \kappa_2\nabla u_2(\xb)\cdot \nabla v_2(\xb)d\xb.
\end{align}
Or equivalently,
\begin{align*}
& NL(u_1,u_2,v_1,v_2,\gamma) - L(u_1,u_2,v_1,v_2)  \\[5mm]
&=  \int\limits_{\Omega_1 \cup\Gamma_1\cup \Gamma_{12}}\int\limits_{\Omega_1\cup\Gamma_1\cup \Gamma_{12}}
    (u_1(\xb)-u_1(\yb))(v_1(\xb)-v_1(\yb)) \gamma_{11}(\xb,\yb)d\yb\,d\xb \\
&+ \int\limits_{\Omega_1\cup\Gamma_1}\int\limits_{\Omg_2}
(u_1(\xb)-u_1(\yb))(v_1(\xb)-v_1(\yb)) \gamma_{12}(\xb,\yb)d\yb\,d\xb \\
&+\int\limits_{\Omg_2 \cup \Gamma_{21}}\int\limits_{\Omg_2 \cup \Gamma_{21}} (u_2(\xb)-u_2(\yb))(v_2(\xb)-v_2(\yb)) \gamma_{22}(\xb,\yb)d\yb\,d\xb \\
&+\int\limits_{\Omg_2}\int\limits_{\Omega_1\cup\Gamma_1} (u_2(\xb)-u_2(\yb))(v_2(\xb)-v_2(\yb)) \gamma_{21}(\xb,\yb)d\yb\,d\xb\\
&- \int_{\Omega_1} \kappa_1\nabla u_1(\xb)\cdot \nabla v_1(\xb)d\xb  
-\int_{\Omega_2} \kappa_2\nabla u_2(\xb)\cdot \nabla v_2(\xb)d\xb \\
&-  \int\limits_{\Omega_1 \cup\Gamma_1}\int\limits_{\Gamma_{12}}
(u_1(\xb)-u_1(\yb))(v_1(\xb)-v_1(\yb)) \gamma_{11}(\xb,\yb)d\yb\,d\xb\\
&-  \int\limits_{\Gamma_{12}}\int\limits_{\Omega_1\cup\Gamma_1\cup \Gamma_{12}}(u_1(\xb)-u_1(\yb))(v_1(\xb)-v_1(\yb)) \gamma_{11}(\xb,\yb)d\yb\,d\xb \\
&-  \int\limits_{\Omega_2}\int\limits_{\Gamma_{21}}
(u_2(\xb)-u_2(\yb))(v_2(\xb)-v_2(\yb)) \gamma_{22}(\xb,\yb)d\yb\,d\xb\\
&- \int\limits_{\Gamma_{21}}\int\limits_{\Omega_2 \cup \Gamma_{21}}
(u_2(\xb)-u_2(\yb))(v_2(\xb)-v_2(\yb)) \gamma_{22}(\xb,\yb)d\yb\,d\xb.
\end{align*}
Note that $\gamma_{ii}(\xb,\yb) = C_{ii}(\delta_i) \chi_{B_i(\mathbf{x})}(\mathbf{y})$ is a radial function, i.e. $\gamma_{ii}(\xb,\yb) = \gamma_{ii}(|\xb-\yb|)$. Let $\mathbf{z} = \mathbf{x}-\mathbf{y}$, then $\gamma_{ii}=\gamma_{ii}(|\mathbf{z}|) = C_{ii}(\delta_i) \chi_{B_i(\mathbf{0})}(\mathbf{z})$. Define $(K^{NL}_i)_{jk} := \int_{B_i(\mathbf{0})}\gamma_{ii}(|\mathbf{z}|)z_j\,z_k d\mathbf{z}$, with $j,k=1,2$.
It follows that
\begin{align}
 (K^{NL}_i)_{11} = \int_{B_i(\mathbf{0})}\gamma_{ii}(|\mathbf{z}|)z_1\,z_1 d\mathbf{z} = \int_{0}^{2\,\pi}\Big(\int_{0}^{\delta_i} C_{ii}(\delta_i) \rho^3 \cos^2(\theta)d \rho \Big)d \theta = 
 C_{ii}(\delta_i)\,\pi\,\dfrac{\delta_i^4}{4} = \kappa_i. 
\end{align}
\begin{align}
 (K^{NL}_i)_{jk} = \int_{B_i(\mathbf{0})}\gamma_{ii}(|\mathbf{z}|)z_j\,z_k d\mathbf{z} = \int_{0}^{2\,\pi}\Big(\int_{0}^{\delta_i} C_{ii}(\delta_i) \rho^3 \cos(\theta)\,\sin(\theta)d \rho \Big)d \theta = 0, \quad j \neq k 
\end{align}
\begin{align}
 (K^{NL}_i)_{22} = \int_{B_i(\mathbf{0})}\gamma_{ii}(|\mathbf{z}|)z_2\,z_2 d\mathbf{z} = \int_{0}^{2\,\pi}\Big(\int_{0}^{\delta_i} C_{ii}(\delta_i) \rho^3 \sin^2(\theta)d \rho \Big)d \theta = 
 C_{ii}(\delta_i)\,\pi\,\dfrac{\delta_i^4}{4} = \kappa_i. 
\end{align}
 Define 
\begin{align*}
 K^L_i := \lim_{\delta_i \rightarrow 0} K_i^{NL}, 
\quad {\rm then} \quad
K^L_i = 
\left[ \begin{matrix}
\kappa_i &   0\\
0 &   \kappa_i
\end{matrix}\right].
\end{align*}
It follows from a result in   \cite{du2012analysis} that
\begin{align*}
\int\limits_{\Omega_1 \cup\Gamma_1\cup \Gamma_{12}}\int\limits_{\Omega_1\cup\Gamma_1\cup \Gamma_{12}} \!\!\!\!\!(u_1(\xb)-u_1(\yb))(v_1(\xb)-v_1(\yb)) \gamma_{11}(\xb,\yb)d\yb\,d\xb = \int_{\Omega_1} \kappa_1\nabla u_1(\xb)\cdot \nabla v_1(\xb)d\xb 
+\mathcal O(\delta_1^2).
\end{align*}
\begin{align*}
\int\limits_{\Omg_2 \cup  \Gamma_{21}}\int\limits_{\Omg_2 \cup \Gamma_{21}} \!\!\!\!\!(u_2(\xb)-u_2(\yb))(v_2(\xb)-v_2(\yb)) \gamma_{22}(\xb,\yb)d\yb\,d\xb = \int_{\Omega_2} \kappa_2\nabla u_2(\xb)\cdot \nabla v_2(\xb)d\xb
+\mathcal O(\delta_2^2).
\end{align*}
Let's now focus on each of the other contributions in \eqref{eq:terms}, one at the time.
\begin{equation}
 \begin{aligned}
  &\int\limits_{\Omega_1\cup\Gamma_1}\int\limits_{\Omg_2}
(u_1(\xb)-u_1(\yb))(v_1(\xb)-v_1(\yb)) \gamma_{12}(\xb,\yb)d\yb\,d\xb \\
&\qquad \qquad= C_{12}(\delta_1)
\int\limits_{\Omega_1\cup\Gamma_1}\int\limits_{\Omg_2 \cap B_1{(\xb)}} (u_1(\xb)-u_1(\yb))(v_1(\xb)-v_1(\yb))  d\yb\,d\xb\nonumber\\
&\qquad \qquad=C_{12}(\delta_1)
\int\limits_{\underline{\Gamma}_{21}}\int\limits_{\Gamma_{12} \cap B_1{(\xb)}} (u_1(\xb)-u_1(\yb))(v_1(\xb)-v_1(\yb))  d\yb\,d\xb \\
&\qquad \qquad \approx 
C_{12}(\delta_1)
\int\limits_{\underline{\Gamma}_{21}}\int\limits_{\Gamma_{12} \cap B_1{(\xb)}} \mathcal{O}(\delta_1)\mathcal{O}(\delta_1) d\yb\,d\xb \quad \mbox{(using Taylor expansions)}\\
&\qquad \qquad\approx C_{12}(\delta_1)\mathcal{O}(\delta_1^2) \int\limits_{\underline{\Gamma}_{21}}\mathcal{O}(\delta_1^2)d\yb \approx  C_{12}(\delta_1)\mathcal{O}(\delta_1^2) \mathcal{O}(\delta_1^2) \mathcal{O}(\delta_1) =  C_{12}(\delta_1) \mathcal{O}(\delta_1^5) = \mathcal{O}(\delta_1).
 \end{aligned}
\end{equation}
In a similar manner, we obtain
\begin{equation}
 \begin{aligned}
  &\int\limits_{\Omega_2}\int\limits_{\Omg_1\cup\Gamma_1}
(u_2(\xb)-u_2(\yb))(v_2(\xb)-v_2(\yb)) \gamma_{21}(\xb,\yb)d\yb\,d\xb\\
&\qquad \qquad=C_{21}(\delta_2)
\int\limits_{\Omega_2}\int\limits_{\Omg_1\cup\Gamma_1 \cap B_2{(\xb)}} (u_2(\xb)-u_2(\yb))(v_2(\xb)-v_2(\yb))  d\yb\,d\xb\nonumber\\
&\qquad \qquad=C_{21}(\delta_2)
\int\limits_{\underline{\Gamma}_{12}}\int\limits_{\Gamma_{21} \cap B_2{(\xb)}} (u_2(\xb)-u_2(\yb))(v_2(\xb)-v_2(\yb))  d\yb\,d\xb\\
&\qquad \qquad\approx C_{21}(\delta_2)\mathcal{O}(\delta_2^2) \int\limits_{\underline{\Gamma}_{12}}\mathcal{O}(\delta_2^2)d\yb \approx  
\mathcal{O}(\delta_2).
 \end{aligned}
\end{equation}
\begin{equation}
 \begin{aligned}
  &\int\limits_{\Omega_1\cup\Gamma_1}\int\limits_{\Gamma_{12}}
(u_1(\xb)-u_1(\yb))(v_1(\xb)-v_1(\yb)) \gamma_{11}(\xb,\yb)d\yb\,d\xb =\\
&\qquad \qquad=C_{11}(\delta_1)
\int\limits_{\Omega_1\cup\Gamma_1}\int\limits_{\Gamma_{12} \cap B_1{(\xb)}} (u_1(\xb)-u_1(\yb))(v_1(\xb)-v_1(\yb))  d\yb\,d\xb\nonumber\\
&\qquad \qquad=C_{11}(\delta_1)
\int\limits_{\underline{\Gamma}_{21}}\int\limits_{\Gamma_{12} \cap B_1{(\xb)}} (u_1(\xb)-u_1(\yb))(v_1(\xb)-v_1(\yb))  d\yb\,d\xb \\
&\qquad \qquad\approx C_{11}(\delta_1)\mathcal{O}(\delta_1^2) \int\limits_{\underline{\Gamma}_{21}}\mathcal{O}(\delta_1^2)d\yb \approx  
\mathcal{O}(\delta_1).
\end{aligned}
\end{equation}
\begin{equation}
\begin{aligned}
&\int\limits_{\Gamma_{12}}\int\limits_{\Omega_1\cup\Gamma_1\cup \Gamma_{12}}
(u_1(\xb)-u_1(\yb))(v_1(\xb)-v_1(\yb)) \gamma_{11}(\xb,\yb)d\yb\,d\xb \\ &\qquad \qquad =C_{11}(\delta_1)
\int\limits_{\Gamma_{12}}\int\limits_{(\Omega_1\cup\Gamma_1\cup \Gamma_{12}) \cap B_1{(\xb)}} \!\!\!\!(u_1(\xb)-u_1(\yb))(v_1(\xb)-v_1(\yb))  d\yb\,d\xb\nonumber\\
&\qquad \qquad=C_{11}(\delta_1)
\int\limits_{\Gamma_{12}}\int\limits_{( \underline{\Gamma}_{21} \cup \Gamma_{12}) \cap B_1{(\xb)}} (u_1(\xb)-u_1(\yb))(v_1(\xb)-v_1(\yb))  d\yb\,d\xb\\
&\qquad \qquad\approx C_{11}(\delta_1)\mathcal{O}(\delta_1^2) \int\limits_{\Gamma_{12}}\mathcal{O}(\delta_1^2)d\yb \approx
\mathcal{O}(\delta_1).
 \end{aligned}
\end{equation}
\begin{equation}
\begin{aligned}
&\int\limits_{\Omega_2}\int\limits_{\Gamma_{21}}
    (u_2(\xb)-u_2(\yb))(v_2(\xb)-v_2(\yb)) \gamma_{22}(\xb,\yb)d\yb\,d\xb \\
    &\qquad \qquad=C_{22}(\delta_2)
\int\limits_{\Omega_2}\int\limits_{\Gamma_{21} \cap B_2{(\xb)}} (u_2(\xb)-u_2(\yb))(v_2(\xb)-v_2(\yb))  d\yb\,d\xb\nonumber\\
&\qquad \qquad=C_{22}(\delta_2)
\int\limits_{\underline{\Gamma}_{12}}\int\limits_{\Gamma_{21} \cap B_2{(\xb)}} (u_2(\xb)-u_2(\yb))(v_2(\xb)-v_2(\yb))  d\yb\,d\xb \\
&\qquad \qquad\approx C_{22}(\delta_2)\mathcal{O}(\delta_2^2) \int\limits_{\underline{\Gamma}_{12}}\mathcal{O}(\delta_2^2)d\yb \approx  
\mathcal{O}(\delta_2).
 \end{aligned}
\end{equation}
\begin{equation}
 \begin{aligned}
  &\int\limits_{\Gamma_{21}}\int\limits_{\Omega_2 \cup \Gamma_{21}}
    (u_2(\xb)-u_2(\yb))(v_2(\xb)-v_2(\yb)) \gamma_{22}(\xb,\yb)d\yb\,d\xb\\
    &\qquad \qquad= C_{22}(\delta_2)
\int\limits_{\Gamma_{21}}\int\limits_{(\Omega_2 \cup \Gamma_{21}) \cap B_2{(\xb)}} (u_2(\xb)-u_2(\yb))(v_2(\xb)-v_2(\yb))  d\yb\,d\xb\nonumber\\
&\qquad \qquad=C_{22}(\delta_2)
\int\limits_{\Gamma_{21}}\int\limits_{(\underline{\Gamma}_{12}\cup\Gamma_{21}) \cap B_2{(\xb)}} (u_2(\xb)-u_2(\yb))(v_2(\xb)-v_2(\yb))  d\yb\,d\xb\\
&\qquad \qquad\approx C_{22}(\delta_2)\mathcal{O}(\delta_2^2) \int\limits_{\Gamma_{21}}\mathcal{O}(\delta_2^2)d\yb \approx
\mathcal{O}(\delta_2).
 \end{aligned}
\end{equation}
Therefore, we conclude that
$
|NL(u_1,u_2,v_1,v_2,\gamma)-L(u_1,u_2,v_1,v_2)| \approx
\mathcal O(\delta_1) + \mathcal O(\delta_2).$


\section*{References}
\begin{thebibliography}{10}
\bibitem{benson2000application}
D.~A. Benson, S.~W. Wheatcraft, M.~M. Meerschaert, Application of a fractional
  advection-dispersion equation, Water resources research 36~(6) (2000)
  1403--1412.

\bibitem{schumer2001eulerian}
R.~Schumer, D.~A. Benson, M.~M. Meerschaert, S.~W. Wheatcraft, Eulerian
  derivation of the fractional advection--dispersion equation, Journal of
  contaminant hydrology 48~(1-2) (2001) 69--88.

\bibitem{schumer2003multiscaling}
R.~Schumer, D.~A. Benson, M.~M. Meerschaert, B.~Baeumer, Multiscaling
  fractional advection-dispersion equations and their solutions, Water
  Resources Research 39~(1).

\bibitem{delgoshaie2015non}
A.~H. Delgoshaie, D.~W. Meyer, P.~Jenny, H.~A. Tchelepi, Non-local formulation
  for multiscale flow in porous media, Journal of Hydrology 531 (2015)
  649--654.

\bibitem{silling2000reformulation}
S.~A. Silling, Reformulation of elasticity theory for discontinuities and
  long-range forces, Journal of the Mechanics and Physics of Solids 48~(1)
  (2000) 175--209.

\bibitem{ha2011characteristics}
Y.~D. Ha, F.~Bobaru, Characteristics of dynamic brittle fracture captured with
  peridynamics, Engineering Fracture Mechanics 78~(6) (2011) 1156--1168.

\bibitem{littlewood2010simulation}
D.~J. Littlewood, Simulation of dynamic fracture using peridynamics, finite
  element modeling, and contact, in: Proceedings of the ASME 2010 International
  Mechanical Engineering Congress and Exposition (IMECE), 2010, pp. 209--217.

\bibitem{buades2010image}
A.~Buades, B.~Coll, J.-M. Morel, Image denoising methods. a new nonlocal
  principle, SIAM review 52~(1) (2010) 113--147.

\bibitem{gilboa2007nonlocal}
G.~Gilboa, S.~Osher, Nonlocal linear image regularization and supervised
  segmentation, Multiscale Modeling \& Simulation 6~(2) (2007) 595--630.

\bibitem{gilboa2008nonlocal}
G.~Gilboa, S.~Osher, Nonlocal operators with applications to image processing,
  Multiscale Modeling \& Simulation 7~(3) (2008) 1005--1028.

\bibitem{lou2010image}
Y.~Lou, X.~Zhang, S.~Osher, A.~Bertozzi, Image recovery via nonlocal operators,
  Journal of Scientific Computing 42~(2) (2010) 185--197.

\bibitem{schekochihin2008mhd}
A.~A. Schekochihin, S.~C. Cowley, T.~A. Yousef, {MHD} turbulence: Nonlocal,
  anisotropic, nonuniversal?, in: IUTAM Symposium on computational physics and
  new perspectives in turbulence, Springer, 2008, pp. 347--354.

\bibitem{askari2008peridynamics}
E.~Askari, F.~Bobaru, R.~Lehoucq, M.~Parks, S.~Silling, O.~Weckner,
  Peridynamics for multiscale materials modeling, in: Journal of Physics:
  Conference Series, Vol. 125, IOP Publishing, 2008, p. 012078.

\bibitem{alali2012multiscale}
B.~Alali, R.~Lipton, Multiscale dynamics of heterogeneous media in the
  peridynamic formulation, Journal of Elasticity 106~(1) (2012) 71--103.

\bibitem{bates1999integrodifferential}
P.~W. Bates, A.~Chmaj, An integrodifferential model for phase transitions:
  stationary solutions in higher space dimensions, Journal of statistical
  physics 95~(5-6) (1999) 1119--1139.

\bibitem{fife2003some}
P.~Fife, Some nonclassical trends in parabolic and parabolic-like evolutions,
  in: Trends in nonlinear analysis, Springer, 2003, pp. 153--191.

\bibitem{meerschaert2011stochastic}
M.~M. Meerschaert, A.~Sikorskii, Stochastic models for fractional calculus,
  Vol.~43, Walter de Gruyter, 2011.

\bibitem{d2017nonlocal}
M.~D'Elia, Q.~Du, M.~Gunzburger, R.~Lehoucq, Nonlocal convection-diffusion
  problems on bounded domains and finite-range jump processes, Computational
  Methods in Applied Mathematics 17~(4) (2017) 707--722.

\bibitem{burch2014exit}
N.~Burch, M.~D'Elia, R.~B. Lehoucq, The exit-time problem for a {M}arkov jump
  process, The European Physical Journal Special Topics 223~(14) (2014)
  3257--3271.

\bibitem{alali2015peridynamics}
B.~Alali, M.~Gunzburger, Peridynamics and material interfaces, Journal of
  Elasticity 120~(2) (2015) 225--248.

\bibitem{seleson2013interface}
P.~Seleson, M.~Gunzburger, M.~L. Parks, Interface problems in nonlocal
  diffusion and sharp transitions between local and nonlocal domains, Computer
  Methods in Applied Mechanics and Engineering 266 (2013) 185--204.

\bibitem{Hansbo_02_CMAME}
A.~Hansbo, P.~Hansbo,
  \href{http://www.sciencedirect.com/science/article/pii/S0045782502005248}{An
  unfitted finite element method, based on nitsche's method, for elliptic
  interface problems}, Computer Methods in Applied Mechanics and Engineering
  191~(47-48) (2002) 5537 -- 5552.
\newblock \href {http://dx.doi.org/10.1016/S0045-7825(02)00524-8}
  {\path{doi:10.1016/S0045-7825(02)00524-8}}.

\bibitem{Chen_11_SINUM}
W.~Chen, M.~Gunzburger, F.~Hua, X.~Wang,
  \href{http://dx.doi.org/10.1137/080740556}{A parallel robin--robin domain
  decomposition method for the stokes--darcy system}, SIAM Journal on Numerical
  Analysis 49~(3) (2011) 1064--1084.
\newblock \href {http://arxiv.org/abs/http://dx.doi.org/10.1137/080740556}
  {\path{arXiv:http://dx.doi.org/10.1137/080740556}}, \href
  {http://dx.doi.org/10.1137/080740556} {\path{doi:10.1137/080740556}}.

\bibitem{Cao_10_CMS}
Y.~Cao, M.~Gunzburger, F.~Hua, X.~Wang,
  \href{http://projecteuclid.org/euclid.cms/1266935011}{Coupled stokes-darcy
  model with beavers-joseph interface boundary condition}, Commun. Math. Sci.
  8~(1) (2010) 1--25.

\bibitem{Javili_14_CMAME}
A.~Javili, S.~Kaessmair, P.~Steinmann, General imperfect interfaces, Computer
  Methods in Applied Mechanics and Engineering 275 (2014) 76 -- 97.

\bibitem{Ming_08_JAP}
M.~Hu, P.~Keblinski, J.-S. Wang, N.~Raravikar,
  \href{https://doi.org/10.1063/1.3000441}{Interfacial thermal conductance
  between silicon and a vertical carbon nanotube}, Journal of Applied Physics
  104~(8) (2008) 083503.
\newblock \href {http://arxiv.org/abs/https://doi.org/10.1063/1.3000441}
  {\path{arXiv:https://doi.org/10.1063/1.3000441}}, \href
  {http://dx.doi.org/10.1063/1.3000441} {\path{doi:10.1063/1.3000441}}.

\bibitem{Berber_00_PRL}
S.~Berber, Y.-K. Kwon, D.~Tom{\'a}nek,
  \href{https://www.scopus.com/inward/record.uri?eid=2-s2.0-0000765076&doi=10.1103%2fPhysRevLett.84.4613&partnerID=40&md5=8e07de64a7a0589a24f0f478e874b5e3}{Unusually
  high thermal conductivity of carbon nanotubes}, Physical Review Letters
  84~(20) (2000) 4613--4616, cited By 2512.
\newblock \href {http://dx.doi.org/10.1103/PhysRevLett.84.4613}
  {\path{doi:10.1103/PhysRevLett.84.4613}}.

\bibitem{Lemarie_15_PCS}
F.~Lemari{\'e}, E.~Blayo, L.~Debreu,
  \href{http://www.sciencedirect.com/science/article/pii/S1877050915012818}{Analysis
  of ocean-atmosphere coupling algorithms: Consistency and stability}, Procedia
  Computer Science 51 (2015) 2066 -- 2075, international Conference On
  Computational Science, \{ICCS\} 2015Computational Science at the Gates of
  Nature.
\newblock \href
  {http://dx.doi.org/http://dx.doi.org/10.1016/j.procs.2015.05.473}
  {\path{doi:http://dx.doi.org/10.1016/j.procs.2015.05.473}}.

\bibitem{Lions_95_JMPA}
J.~L. Lions, R.~Temam, S.~Wang, Mathematical theory for the coupled
  atmosphere-ocean models (cao iii), Journal de Mathematiques Pures et
  Appliquees 74~(2) (1995) 105--163.

\bibitem{Kwak_17_IJNME}
D.~Kwak, S.~Lee, Y.~Hyon, A new finite element for interface problems having
  robin type jump, Int. J. Num. Analysis and Modeling 14~(4-5) (2017) 532--549.

\bibitem{Du_13_MMMAS}
Q.~Du, M.~Gunzburger, R.~B. Lehoucq, K.~Zhou, A nonlocal vector calculus,
  nonlocal volume-constrained problems, and nonlocal balance laws, Mathematical
  Models and Methods in Applied Sciences 23~(03) (2013) 493--540.

\bibitem{Mahan_09_PRB}
G.~D. Mahan, \href{https://link.aps.org/doi/10.1103/PhysRevB.79.075408}{Kapitza
  thermal resistance between a metal and a nonmetal}, Phys. Rev. B 79 (2009)
  075408.
\newblock \href {http://dx.doi.org/10.1103/PhysRevB.79.075408}
  {\path{doi:10.1103/PhysRevB.79.075408}}.

\bibitem{du2012analysis}
Q.~Du, M.~Gunzburger, R.~B. Lehoucq, K.~Zhou, Analysis and approximation of
  nonlocal diffusion problems with volume constraints, SIAM review 54~(4)
  (2012) 667--696.
\end{thebibliography}
\end{document}